\documentclass[11pt]{amsart}
\usepackage{amscd,amssymb}
\usepackage{amsfonts}
\usepackage{amsmath}
\usepackage{graphicx} 
\usepackage{color}
\usepackage{hyperref} 

\newcommand{\mrm}{\mathrm}
\newcommand{\mbb}{\mathbb}
\newcommand{\mf}{\mathfrak}

\newcommand{\Uq}{U_q(\mathfrak{sl}_2)}
\newcommand{\D}{\Delta}
\newcommand{\ot}{\otimes}
\newcommand{\ld}{\lambda}
\newcommand{\al}{\alpha}
\newcommand{\Zp}{\mathbb{Z}_{\geq 0}}
\newcommand\q[1]{\left[ #1 \right]_q}
\newcommand\qt[6]{
	\begin{bmatrix}
		#1 & #2 & #3\\
		#4 & #5 & #6
	\end{bmatrix}_q
}
\newcommand\qtm[6]{
	\begin{bmatrix}
		#1 & #2 & #3\\
		#4 & #5 & #6
	\end{bmatrix}_{q^{-1}}
}
\newcommand\qs[6]{
	\begin{Bmatrix}
		#1 & #2 & #3\\
		#4 & #5 & #6
	\end{Bmatrix}_q
}
\newcommand{\adm}{\mrm{ADM}}
\newcommand\p[2]{\prod\limits_{j=0}^{#1}\q{#2}}
\newcommand\qbin[2]{\left[  \begin{array}{c}  #1 \\ #2  \end{array}  \right]_q}
\newcommand\dd[1]{\Delta(#1)}

\newtheorem{theorem}{Theorem}[section]

\newtheorem{lemma}[theorem]{Lemma}
\newtheorem{definition}[theorem]{Definition} 
\newtheorem{conjecture}[theorem]{Conjecture}

\newtheorem{remark}[theorem]{Remark}
\title[Quantum $6j$-symbols for Verma modules]{Infinite-dimensional representations of $\Uq$ and the shadow world: quantum $6j$-symbols for Verma modules}
\author{Dmitry Solovyev}
\address{D.S.: Yau Mathematical Sciences Center, Tsinghua University, Beijing, China}
\email{dimsol42@gmail.com}

\begin{document}
\begin{abstract}
This paper initiates the study of invariants of links associated to infinite-dimensional representations of $\Uq$ using graphical representation for quantum $6j$-symbols, the shadow world. We obtain formulae for $q3j$-symbols and $q6j$-symbols for Verma modules and study their properties. This hints at the structure of a possible target category for the Reshetikhin-Turaev functor for such an invariant.
\end{abstract}
\maketitle
\tableofcontents

\addtocontents{toc}{\protect\setcounter{tocdepth}{1}}
\section{Introduction}
\subsection{Background} Quantum $6j$-symbols are fundamental objects in the representation theory of quantum groups, with deep connections to diverse fields of research such as $3$-manifold invariants \cite{TV,GKT}, $q$-orthogonal polynomials \cite{AW,GR}, conformal field theory \cite{MS,Te,TeVa}, and models of $3$d quantum gravity \cite{PR,CEZ}. In the present work, we focus on their role in knot theory, where they provide the algebraic framework underlying quantum knot invariants. In \cite{KR} Kirillov and Reshetikhin proposed a functional based on graphical representation for $q6j$-symbols of finite-dimensional irreducible representations of $\Uq$, also known as the shadow world. This functional is inspired by the solid-on-solid model from statistical mechanics \cite{JMO} and allows one to compute the colored Jones polynomial. It is equivalent to the functional constructed via the $\mathcal{R}$-matrix approach \cite{R1}, which was later extended to a functorial framework known as the Reshetikhin-Turaev functor \cite{RT}. Building on this, Reshetikhin and Turaev \cite{RT2} constructed a non-perturbative completion of the Chern-Simons partition function: a surgery topological quantum field theory that upgrades the link invariants defined by the RT functor to invariants of $3$-manifolds, now known as the Witten-Reshetikhin-Turaev invariants.\par

Works \cite{GPV, GPPV} predict new 3-manifold invariants $\hat{Z}_a(M_3;q)$, defined via the 3d $\mathcal{N}=2$ theory $T[M_3]$ obtained by compactifying $6$d $(2,0)$ M-theory on $M_3$. These invariants take the form of $q$-series convergent in the unit disk, whose radial limits at roots of unity recover the WRT invariants. Gukov and Manolescu \cite{GM} further conjectured that for each knot $K$ there exists a two-variable series $F_K(x,q)$, obtained by resumming the Melvin-Morton-Rozansky expansion \cite{MM,Roz} of the colored Jones polynomial, with surgery formulas relating $F_K(x,q)$ to $\hat{Z}_a(M_3;q)$. These developments suggest that $\hat{Z}_a(M_3;q)$ admits a non-perturbative completion in the form of a surgery TQFT, constructed from the RT functor with a suitable target category. A key step in this direction is a functorial formulation of $F_K(x,q)$, which requires identifying the appropriate categorical target for the RT functor. This would also pave the way for a categorification analogous to Khovanov homology \cite{Kh} and further solidify the deep connections of these invariants with logarithmic vertex algebras \cite{CCF}, holomorphic curve counting \cite{EGG}, geometric representation theory \cite{GHN} and other fields of research.\par

The equivalence between the $\mathcal{R}$-matrix (edge-state) model and the shadow world (face-state) model can be understood via the bulk-boundary correspondence: in the edge-state model, degrees of freedom are localized along strands with interactions determined by intertwiners and braiding, while in the shadow world (face-state) model, degrees of freedom reside in the bulk with interactions governed by fusion rules. The correspondence asserts that edge scattering data can be recast as fusion constraints on faces, providing dual descriptions of the same representation-theoretic structure. We use this approach to gain insights into the Grothendieck ring of a potential categorical target for the RT functor that recovers $F_K(x,q)$.\par

The studies \cite{P1,P2} initiated the investigation of $F_K(x,q)$ via the $\mathcal{R}$-matrix approach for Verma modules of $\Uq$. Hence, it is natural to explore possible extensions of the bulk-boundary correspondence beyond finite-dimensional representations of $\Uq$ and initiate the study of knot and link invariants arising from Verma modules using the shadow world. As opposed to the $\mathcal{R}$-matrix approach, which relies on the summation over the weight spaces of representations coloring the strands of a knot, the functional given by the shadow world relies on the summation over the highest weights in the tensor product decomposition of the representations coloring the strands and the gleams. Combined with known results on $F_K(x,q)$, this provides insights into the Grothendieck ring of a potential target category for the RT functor.\par

\subsection{Results} This paper establishes formulae and properties for quantum Clebsch-Gordan coefficients($q3j$-symbols) and quantum Racah-Wigner coefficients($q6j$-symbols) associated with a braided monoidal version of the $q$-analog of BGG category $\mathcal{O}$ \cite{BG, BGG}, which includes Verma modules of the highest weight \cite{DCRS}. Our focus is on tensor products of Verma modules $M_\ld$, restricting attention to generic combinations of highest weights to ensure irreducibility of the components in the tensor product decomposition. Guided by the bulk-boundary correspondence, we discuss the equivalence of the $\mathcal{R}$-matrix and the shadow world approaches for Verma modules. We introduce a functional $\mathsf{SW}_\infty$ and investigate its potential as a model for the Gukov-Manolescu series $F_K(x,q)$. Based on this discussion, we formulate a conjecture for a target category for the RT functor and its fusion rules.\par

In Section \ref{sectionUq} we introduce necessary definitions and notions related to Verma modules of $\Uq$ for generic values of $q$. We focus on the case where the highest weights $\ld_j \in \mbb{C}$ satisfy $\sum_j\ld_j\notin\frac{1}{2}\mbb{Z}_{\geq 0}$, ensuring all modules are irreducible and
\begin{equation*}
M_{\ld_1}\ot M_{\ld_2}=\bigoplus\limits_{J=0}^\infty M_{\ld_1+\ld_2-J}.
\end{equation*}\par

In Section \ref{sectionq3j} we derive formulae for $q3j$-symbols. Embeddings $\psi_\ld: M_\ld\hookrightarrow M_{\ld_1}\ot M_{\ld_2}$ are defined by coefficients given by
\begin{eqnarray*}
    \qt{\ld_1}{\ld_2}{\ld}{a_1}{a_2}{a}^\psi=(-1)^{\frac{J}{2}+\ld_1-a_1}q^{(\ld_2^2-a_2^2)-(\ld^2-a^2)-(\ld_1-a_1)}\frac{\q{\ld-a}!}{\prod\limits_{j=1}^{\ld-a}\q{\ld+a+j}}\times\\
    \times \dd{\ld_1,\ld_2,\ld}\sum\limits_{z,w}(-1)^zq^{z(\ld+a+1)}\frac{\prod\limits_{j=1}^{z}\q{\ld_1+a_1+j}\prod\limits_{j=1}^{w}\q{\ld_2+a_2+j}}{\q{z}!\q{w}!\q{\ld_1-a_1-z}! \q{\ld_2-a_2-w}!},\nonumber
\end{eqnarray*}
where $(z,w)\in\{(0,\ld-a),(1,\ld-a-1),\ldots,(\ld-a,0)\}$. In the formula above 
$$
\dd{\ld_1,\ld_2,\ld}=\left(\q{\ld_1+\ld_2-\ld}!\prod\limits_{j=0}^{\ld_1+\ld_2-\ld-1}\frac{\q{2\ld_1-j}\q{2\ld_2-j}}{\q{2\ld+2+j}}\right)^{\frac{1}{2}}
$$
and $a=\ld-k$, $a_j=\ld_j-k_j$, for some $k,k_j\in\mbb{Z}_{\geq 0}$, $j=1,2$. This is an infinite-dimensional analog of the Racah-Fock formula. In Section \ref{qtprop}, we use string diagrams, following \cite{KR}, to establish properties of $q3j$-symbols, including orthogonality relations and symmetries.\par

In Section \ref{sectionq6j} we use results obtained for $q3j$-symbols to derive formulae for $q6j$-symbols, with particular emphasis on the implementation of the Racah method \cite{BL}. Coefficients in the associator for the tensor product decomposition of $M_{\ld_1}\ot M_{\ld_2}\ot M_{\ld_3}$ are given by
\begin{eqnarray*}
    &\qs{\ld_1}{\ld_2}{\ld_{12}}{\ld_3}{\ld}{\ld_{23}}= \frac{\dd{\ld_1,\ld_{23},\ld}}{\prod\limits_{j=0}^{\ld_1+\ld_{23}-\ld-1}\q{2\ld_{23}-j}}\frac{\dd{\ld_2,\ld_3,\ld_{23}}}{\prod\limits_{j=0}^{\ld_2+\ld_3-\ld_{23}-1}\q{2\ld_3-j}}\frac{\dd{\ld_{12},\ld_3,\ld}}{\prod\limits_{j=0}^{\ld_{12}+\ld_3-\ld-1}\q{2\ld_{12}-j}}\frac{\dd{\ld_1,\ld_2,\ld_{12}}}{\prod\limits_{j=0}^{\ld_1+\ld_2-\ld_{12}-1}\q{2\ld_1-j}}\times\nonumber\\
    &\times \frac{1}{\prod\limits_{j=0}^{\ld_1+\ld_2-\ld_{12}-1}\q{2\ld_2-j}}\sum\limits_{z=0}^{\ld_1+\ld_2-\ld_{12}}(-1)^z \prod\limits_{j=0}^{\ld_1+\ld_2+\ld_3-\ld-z-1}\q{2\ld+2+j}\times\\
    &\times\frac{\prod\limits_{j=0}^{z-1}\q{\ld+\ld_1-\ld_{23}-j}\prod\limits_{j=0}^{z-1}\q{\ld-\ld_{12}+\ld_3-j}\prod\limits_{j=0}^{\ld_1+\ld_2-\ld_{12}-z-1}\q{\ld_2+\ld_{23}-\ld_3-j}}{\q{z}!\q{\ld_1+\ld_2-\ld_{12}-z}!\q{\ld_2+\ld_3-\ld_{23}-z}!\q{\ld_{12}+\ld_{23}-\ld_2-\ld+z}!}\nonumber.
\end{eqnarray*}
In the formula above $\ld=\ld_1+\ld_2+\ld_3-J$, for some $J\in\mbb{Z}_{\geq 0}$, and similarly $\ld_{mn}=\ld_m+\ld_n-J_{mn}$, for some $J_{mn}\in\mbb{Z}_{\geq 0}$, where $m,n=1,2,3$. In Section \ref{sectionpropq6j}, we use the shadow world graphical calculus to establish properties of $q6j$-symbols, including orthogonality relation, Racah, Biedenharn-Elliott, and Yang-Baxter identities.\par

With appropriate restrictions imposed on the admissible values of $J$ and on the weights $a_j$, the resulting formulae for $q3j$- and $q6j$-symbols can be applied to tensor products whose components are finite-dimensional irreducible representations, provided that the highest weights are chosen generically and the decomposition involves only irreducible representations. In the above formulae, all $q$-factorial arguments are integer numbers. We set $\q{0}!:= 1$, and terms with negative integer arguments are omitted. Likewise, for any $\gamma \in \mathbb{C}$ we set $\prod_{j=0}^{-1}\q{\gamma\pm j}:=1$ to be consistent with the initial conditions of the recursion relations. Otherwise, we use the products defined in (\ref{qpf}). Many technical details were moved to the Appendix \ref{appp}. They consist of intermediate identities needed to implement the Racah method, which we used to derive formulae for $q6j$-symbols.\par  

Section \ref{swsec} presents preliminary results, with rigorous treatment deferred to \cite{S}. Motivated by the bulk–boundary correspondence for the finite-dimensional irreducible representations of $\Uq$, we tentatively assume that the $\mathcal{R}$-matrix and the shadow world approaches for Verma modules remain equivalent in a suitable sense. We propose a functional $\mathsf{SW}_\infty$, which is invariant with respect to the topological moves studied in Section \ref{sectionq6j}. While its application to Verma modules encounters certain topological limitations, we formulate a conjecture to address these difficulties. In particular, we conjecture a target category for the RT functor and its fusion rules. This category is realized as a quotient of a category of infinite-dimensional representations of $\widetilde{\Uq}$, the Chevalley extended $\Uq$, where we include the $q$-analog of the Weyl element $w$. We refer to these representations as \textit{polarized modules}, which are essentially a Weyl-invariant analytic continuation of finite-dimensional $\Uq$-modules, which preserves its structural properties. Using this category, we modify $\mathsf{SW}_\infty$ and examine $F_K(x,q)$ through the example of the trefoil knot. We also briefly discuss the functoriality of the constructed functional.\par

In conclusion, we outline a strategy for addressing the proposed conjecture, to be investigated in future work. Building on these results, we also plan to develop a rigorous graphical calculus for the proposed category and establish $\mathsf{SW}_\infty$ as a functor. Further open questions are also discussed at the end of this section.\par

\subsection{Acknowledgements} We are grateful to Anatol Kirillov, Nicolai Reshetikhin, Mikhail Khovanov, Igor Frenkel, Semeon Arthamonov, Sergei Gukov for useful comments and insightful discussions. We also thank Sunghyuk Park, Pedro Guicardi and Mrunmay Jagadale for inspiring conversations. The author was supported by the Xing Hua Scholarship program of Tsinghua University.

\newpage


\section{Algebra $\Uq$ and its representations}\label{sectionUq}
\subsection{Hopf algebra $\Uq$}
Given classical Lie algebra $\mathfrak{sl}_2$ generated by $E$, $F$, $H$ with relations
\begin{equation*}
    [H,E]=2E,\quad [H,F]=-2F,\quad [E,F]=H,
\end{equation*}
we consider the quantized universal enveloping algebra $\Uq$, where $q$ is generic. It is generated by elements $E$, $F$, $K^{\pm 1}$ subject to the following relations
$$
KE=q^2 EK,\quad KF=q^{-2}FK,\quad KK^{-1}=K^{-1}K=1,
$$
$$
[E,F]=\frac{K-K^{-1}}{q-q^{-1}}.
$$
This is a Hopf algebra with coproduct $\D$, counit $\epsilon$ and antipode $S$ given by
$$
\D E= E\ot K+ 1\ot E,\quad \D F= F\ot 1 +K^{-1}\ot F,\quad \D K^{\pm1}=K^{\pm 1}\ot K^{\pm 1},
$$
$$
\epsilon(E)=0,\quad \epsilon(F)=0,\quad \epsilon(K^{\pm 1})=1,
$$
$$
S(E)=-EK^{-1},\quad S(F)=-KF,\quad S(K^{\pm 1})=K^{\mp 1}.
$$
\subsection{The universal $R$-matrix}
By $\mathcal{R}\in\Uq^{\otimes 2}$ we denote the universal $R$-matrix. It satisfies properties
\begin{equation}\label{Rid1}
    (\Delta\otimes\mrm{Id}_{\Uq})\mathcal{R}=\mathcal{R}_{13}\mathcal{R}_{23},
\end{equation}
\begin{equation}\label{Rid2}
    (\mrm{Id}_{\Uq}\otimes \Delta)\mathcal{R}=\mathcal{R}_{13}\mathcal{R}_{12},
\end{equation}
\begin{equation}\label{antipodeR}
    (S\otimes \mrm{Id}_{\Uq})\mathcal{R}=\mathcal{R}^{-1},
\end{equation}
where indices show the embeddings of $\mathcal{R}$ into $\Uq^{\otimes 3}$. Identities (\ref{Rid1}) or (\ref{Rid2}) imply Yang-Baxter identity
\begin{equation}\label{YB}
    \mathcal{R}_{12}\mathcal{R}_{13}\mathcal{R}_{23}=\mathcal{R}_{23}\mathcal{R}_{13}\mathcal{R}_{12}.
\end{equation}
Universal $R$-matrix is given by
\begin{equation}\label{rmatr}
    \mathcal{R}=\mrm{exp}\left(\frac{h}{4}H\otimes H\right)\sum\limits_{n\geq 0}\frac{(q-q^{-1})^n}{\q{n}!}q^{\frac{n(n-1)}{2}}E^n\otimes F^n,
\end{equation}
where $q=e^{h/2}$. For the definition of the quantum number $\q{n}$ we refer to Appendix \ref{qnumb}. In what follows we will be concerned with Verma modules, which implies that we have to consider a completion of $\Uq^{\otimes 2}$ suitable for tensor products in the $q$-analog of the BGG category $\mathcal{O}$. For exposition to this subject we refer to Section 1.4 of \cite{DCRS}.

\subsection{Representations of $\Uq$}
Let $\mf{b}$ be the Borel subalgebra of $\mf{sl}_2$ generated by $E$ and $H$, and $\mf{h}$ a Cartan subalgebra generated by $H$. For $\lambda\in\mf{h}^*$ consider one-dimensional $U_q(\mf{b})$-module $\mbb{C}_\ld=\mbb{C} 1_\ld$ defined by $E 1_\ld=0$, $K 1_\ld=q^{\ld(H)} 1_\ld$. The Verma module $M_\ld$ of the highest weight $\ld$ is defined as
$$
M_\ld=\mrm{Ind}_{U_q(\mf{b})}^{\Uq} \mbb{C}_\ld.
$$
In what follows we will abuse notation and instead of $\ld(H)$ simply write $2\ld$ for an arbitrary complex number $\ld$. When $\ld\notin\frac{1}{2}\mbb{Z}_{\geq 0}$, we call such $\ld$ {\it generic}.
\begin{theorem}
Module $M_\ld$ is irreducible iff $\ld$ is generic.
\end{theorem}
When $\ld$ is not generic, the quotient representation
\begin{equation}\label{irrep}
V_\ld=M_\ld/ M_{-\ld-1}
\end{equation}
is the $(2\ld+1)$-dimensional irreducible $\Uq$-module.

Consider linear algebra anti-involution $\rho:\Uq\to \Uq^{\mrm{op}}$ defined as
\begin{equation}\label{rho}
\rho(A)=S(\tau(A)),\quad \forall A\in\Uq,
\end{equation}
where
$$
\tau(E)=F,\quad \tau(F)=E,\quad \tau(K^{\pm 1})=K^{\mp 1}.
$$
It is easy to see that 
$$
\Delta(\rho(A))=(\rho\ot \rho)\Delta A,\quad\forall A\in\Uq,
$$
so this is an automorphism of the coalgebra structure. Let $\pi_{M_\ld}$ be a representation map for Verma module $M_\ld$, and $\{e^\ld_a\}$, where $a=\ld-k$, $k\in\mbb{Z}_{\geq 0}$ be the weight basis.
\begin{definition}\label{defshap}
The Shapovalov form
$$
(-,-):M_\ld\times M_\ld \to \mbb{C}(q)
$$
is a unique bilinear form such that
\begin{itemize}
    \item $(e^\ld_\ld,e^\ld_\ld)=1$
    \item $(\pi(A)e^\ld_a,e^\ld_b)=(e^\ld_a,\pi(\rho(A))e^\ld_b)$, $\quad\forall A\in\Uq$
    \item $c(e^\ld_a,e^\ld_b)=(ce^\ld_a,e^\ld_b)=(e^\ld_a,ce^\ld_b)$,$\quad\forall c\in\mbb{C}(q)$
\end{itemize}
\end{definition}
Now we define the dual Verma module $M_\ld^*$. It is of the highest weight $\lambda$, since we used involution $\tau$ in (\ref{rho}), but one should note that the Shapovalov form is no longer $\Uq$-invariant and, therefore, is not a morphism in the $q$-analog of the category $\mathcal{O}$. Let $\pi_{M_\ld^*}$ be the representation map and $\{f^\ld_a\}$, where $a=\ld-k$, $k\in\mbb{Z}_{\geq 0}$, be the weight basis of $M_\ld^*$. For any $x\in M_\ld$ the Shapovalov form defines a map $\phi:M_\ld\to M_\ld^*$ by the following rule
\begin{equation}\label{isodual}
\phi(x)(-)=(x,-).
\end{equation}
Then
\begin{equation}\label{dualmap}
(\pi_{M_\ld^*}(A)f^\ld_a)(e^\ld_b)=f^\ld_a(\pi_{M_\ld}(\rho(A))e^\ld_b).
\end{equation}
For $\ld$ not generic, map $\phi$ is a homomorphism given by the composition
\begin{equation}\label{seqq}
M_\ld \twoheadrightarrow V_\ld \to V_\ld^*\hookrightarrow M_\ld^*,
\end{equation}
where the dual irreducible $\Uq$-module $V_\ld^*$ can be defined in a similar fashion. Throughout what follows, we will be focusing on Verma modules $M_\ld$, $M_\ld^*$ for generic values of $\ld$.

\subsection{Representation maps $\pi_{M_\ld}$ and $\pi_{M^*_\ld}$}
\subsubsection{Representation map $\pi_{M_\ld}$}
For a given $\ld\in\mbb{C}$ the representation map $\pi_{M_\ld}$ realizing $M_\ld$ is given by
\begin{equation}\label{map}
\begin{cases}
\pi_{M_\ld}(K)e^\ld_a=q^{2a} e^\ld_a\\
\pi_{M_\ld}(E)e^\ld_a=\q{\ld-a} e^\ld_{a+1} \\
\pi_{M_\ld}(F)e^\ld_a=\q{\ld+a} e^\ld_{a-1}
\end{cases}
\end{equation}
where $a=\ld-k$, $k\in\mbb{Z}_{\geq 0}$. It will be convenient for us to write
$$
\begin{cases}
\pi_{M_\ld}(K)e^\ld_{\ld-k}=q^{2(\ld-k)} e^\ld_{\ld-k}\\
\pi_{M_\ld}(E)e^\ld_{\ld-k}=\q{k} e^\ld_{\ld-k+1} \\
\pi_{M_\ld}(F)e^\ld_{\ld-k}=\q{2\ld-k} e^\ld_{\ld-k-1}
\end{cases}
$$
From (\ref{irrep}) we see, that when $\ld$ is not generic the representation map $\pi_{V_\ld}$ is realized by identities (\ref{map}), but for weight vectors $\{e^\ld_a\}$, where $a=\ld,\ld-1,\ldots,-\ld$.

\subsubsection{Representation map $\pi_{M^*_\ld}$}
Consider weight vectors $\{f^\ld_{\ld-k}\}$, $k\in\mbb{Z}_{\geq 0}$. We have
\begin{equation}\label{c1}
f^\ld_{\ld-k}(e^\ld_{\ld-k^\prime})=\delta_{k,k^\prime}.
\end{equation}
Since (\ref{rho}) gives
$$
\rho(E)=-KF,\quad\rho(F)=-EK^{-1},\quad\rho(K^{\pm 1})=K^{\pm 1},
$$
from (\ref{dualmap}) we get
$$
(\pi_{M_\ld^*}(K)f^\ld_{\ld-k})(e^\ld_{\ld-k^\prime})=q^{2(\ld-k^\prime)}\delta_{k,k^\prime},
$$
$$
(\pi_{M_\ld^*}(E)f^\ld_{\ld-k})(e^\ld_{\ld-k^\prime})=-q^{2(\ld-k^\prime-1)}\q{2\ld-k^\prime}\delta_{k,k^\prime+1},
$$
$$
(\pi_{M_\ld^*}(F)f^\ld_{\ld-k})(e^\ld_{\ld-k^\prime})=-q^{-2(\ld-k^\prime)}\q{k^\prime}\delta_{k,k^\prime-1}.
$$
It follows that the representation map $\pi_{M^*_\ld}$ is given by
\begin{equation}\label{dualmap1}
\begin{cases}
\pi_{M^*_\ld}(K)f^\ld_{\ld-k}=q^{2(\ld-k)} f^\ld_{\ld-k}\\
\pi_{M^*_\ld}(E)f^\ld_{\ld-k}=-q^{2(\ld-k)}\q{2\ld-k+1} f^\ld_{\ld-k+1} \\
\pi_{M^*_\ld}(F)f^\ld_{\ld-k}=-q^{-2(\ld-k-1)}\q{k+1} f^\ld_{\ld-k-1}
\end{cases},
\end{equation}
and $\pi_{M^*_\ld}(E)f^\ld_{\ld}=0$. It is easy to see that it is a $\Uq$-module, i.e.
$$
\pi_{M^*_\ld}([E,F])=\q{k+1}\q{2\ld-k}-\q{2\ld-k+1}\q{k}=\q{2(\ld-k)}=\pi_{M^*_\ld}\left(\frac{K-K^{-1}}{q-q^{-1}}\right),
$$
where we have omitted writing weight vector $f^\ld_{\ld-k}$. Similarly, other relations are satisfied.

\subsubsection{Map $\phi:M_\ld \to M^*_\ld$}
Now let us derive how map (\ref{isodual}) acts in a given weight basis. We have
\begin{equation}\label{isodual2}
\phi(\pi_{M_\ld}(A)e^\ld_{\ld-k})=\pi_{M^*_\ld}(A)\phi(e^\ld_{\ld-k}),\quad \forall A\in\Uq
\end{equation}
and suppose
$$
\phi(e^\ld_{\ld-k})=\alpha_{\ld, k} f^\ld_{\ld-k}
$$
for some $\alpha_{\ld,k}$. From the action of $E$ and (\ref{isodual2}) we get the following recursion relation for $\alpha_{\ld,k}$
\begin{equation}\label{recmor}
\q{k} \alpha_{\ld, k-1}=-q^{2(\ld-k)}\q{2\ld-k+1}\alpha_{\ld, k}.
\end{equation}
Condition $(e^\ld_\ld,e^\ld_\ld)=1$ in Definition \ref{defshap} together with (\ref{c1}) fixes the normalization $\alpha_{\ld,0}=1$, and we get the following solution to the recursion relation
\begin{equation}\label{alphaa}
\alpha_{\ld,k}=\frac{(-1)^k q^{k(k-2\ld+1)} \q{k}!}{\p{k-1}{2\ld-j}}.
\end{equation}
Note that we have also obtained
$$
(e^\ld_{\ld-k},e^\ld_{\ld-k^\prime})=\alpha_{\ld,k}\delta_{k,k^\prime}.
$$
We can use map $\phi$ and sequence (\ref{seqq}) to define representation map $\pi_{V^*_\ld}$, which, again, has a similar form to (\ref{dualmap1}).

\subsubsection{The universal $R$-matrix in representation $\pi_{M_\ld}$}
From (\ref{rmatr}) and (\ref{map}) it is easy to see that matrix element $\mathcal{R}^{\ld_1,\ld_2}:=(\pi_{M_{\ld_1}}\otimes \pi_{M_{\ld_2}})\mathcal{R}$ is given by
\begin{eqnarray}\label{rmatrel}
    &\left(\mathcal{R}^{\ld_1,\ld_2}\right)^{a_1,a_2}_{a_1+n,a_2-n}=q^{2(a_1+n)(a_2-n)+\frac{n(n-1)}{2}}\frac{(q-q^{-1})^n}{\q{n}!}\prod\limits_{j=0}^{n-1}\q{\ld_1-a_1-j}\prod\limits_{j=0}^{n-1}\q{\ld_2+a_2-j}.
\end{eqnarray}
From (\ref{YB}) the following identity is satisfied
\begin{equation}\label{YBR}
    \mathcal{R}^{\ld_1,\ld_2}_{12}\mathcal{R}^{\ld_1,\ld_2}_{13}\mathcal{R}^{\ld_1,\ld_2}_{23}=\mathcal{R}^{\ld_1,\ld_2}_{23}\mathcal{R}^{\ld_1,\ld_2}_{13}\mathcal{R}^{\ld_1,\ld_2}_{12}.
\end{equation}
From (\ref{antipodeR}) we get
\begin{eqnarray}
    \left(\left(\mathcal{R}^{\ld_1,\ld_2}\right)^{-1}\right)^{a_1,a_2}_{a_1+n,a_2-n}=(-1)^n q^{-2a_1 a_2+\frac{n(n-1)}{2}-n(n-1)}\times\\
    \frac{(q-q^{-1})^n}{\q{n}!}\prod\limits_{j=0}^{n-1}\q{\ld_1-a_1-j}\prod\limits_{j=0}^{n-1}\q{\ld_2+a_2-j}.\nonumber
\end{eqnarray}
Similar formulae can be obtained for representation maps $\pi_{M^*_\ld}$ but in what follows we will not need it. 
\section{Quantum $3j$-symbols for Verma modules}\label{sectionq3j}
We will proceed with derivation of $q3j$-symbols, also known as quantum Clebsch-Gordan coefficients. We will be following the strategy explored in \cite{K}. The results of this section are summarized in Section \ref{qtprop}.
\subsection{Definition of the $q3j$-symbols}
In what follows we will be mainly concerned with tensor products of Verma modules $M_{\ld_j}$ for generic values of $\ld_j$ such that $\sum_{j} \ld_j$ is also generic. In this case we have the following tensor product decomposition
\begin{equation}\label{decomp}
M_{\ld_1}\ot M_{\ld_2}=\bigoplus\limits_{J=0}^\infty M_{\ld_1+\ld_2-J}.
\end{equation}
For $\nu,\mu\in\mbb{C}$ and $n\in\mbb{Z}$ define the set of admissible pairs as
$$
\adm(\nu,\mu,n):=\{(\nu,\mu-n),\ldots,(\nu-n,\mu)\},
$$
and if all three arguments coincide we write
$$
\adm(n):=\adm(n,n,n).
$$
\begin{definition}
    Embeddings $\psi_\ld:M_\ld\hookrightarrow M_{\ld_1}\ot M_{\ld_2}$ in the decomposition (\ref{decomp}) in terms of weight basis $\{e^\ld_a\}$ are given by formulae
    \begin{equation}\label{emb}
    \psi_\ld(e^\ld_a)=\sum\limits_{a_1,a_2}\qt{\ld_1}{\ld_2}{\ld}{a_1}{a_2}{a}^\psi e^{\ld_1}_{a_1}\ot e^{\ld_2}_{a_2},
    \end{equation}
    where summation in the right hand side is over $(a_1,a_2)\in\adm(\ld_1,\ld_2,\ld_1+\ld_2-a)$. Coefficients near tensor monomials $e^{\ld_1}_{a_1}\ot e^{\ld_2}_{a_2}$ are called $q3j$-symbols.
\end{definition}
Similarly, there are $q3j$-symbols for projections $\pi_\ld:M_{\ld_1}\ot M_{\ld_2}\to M_\ld$
\begin{equation}\label{prj}
    \pi_\ld(e^{\ld_1}_{a_1}\ot e^{\ld_2}_{a_2})= \qt{\ld_1}{\ld_2}{\ld}{a_1}{a_2}{a_1+a_2}^\pi e^\ld_{a_1+a_2}.
\end{equation}
Now we need to compute coefficients in formulae (\ref{emb}) and (\ref{prj}) explicitly.

\subsection{Computing $q3j$-symbols}\label{deriveformulaq3j}
In the following subsections we will be obtaining and solving recursion relation for the $q3j$-symbols.
\subsubsection{Computing $q3j$-symbols for the highest weight vectors}
We will start with deriving formula for
$$
\qt{\ld_1}{\ld_2}{\ld}{a_1}{a_2}{\ld}^\psi.
$$
As we know from (\ref{decomp}), $\ld=\ld_1+\ld_2-J$ for some fixed $J\in\mbb{Z}_{\geq 0}$. We also know, that
$$
a=\ld-n,\quad a_1=\ld_1-n_1,\quad a_2=\ld_2-n_2
$$
for some $n,n_1,n_2\in\mbb{Z}_{\geq 0}$, such that $a=a_1+a_2$ or, equivalently, $J+n=n_1+n_2$. So
$$
(a_1,a_2)\in\{(\ld_1,-\ld_1+a),\dots,(-\ld_2+a,\ld_2)\}=\adm(\ld_1,\ld_2,\ld_1+\ld_2-a)
$$
or, equivalently,
$$
(n_1,n_2)\in\adm(J+n).
$$
Since for highest weight vectors $n=0$, we can rewrite (\ref{emb}) as
\begin{equation}\label{hwv}
    e^\ld_\ld=\sum\limits_{k=0}^J\qt{\ld_1}{\ld_2}{\ld}{\ld_1-k}{\ld_2-J+k}{\ld} e^{\ld_1}_{\ld_1-k}\ot e^{\ld_2}_{\ld_2-J+k},
\end{equation}
where for simplicity we have omitted writing embedding map and embedding index $\psi$. In what follows we will also omit writing representation map $\pi_{M_\ld}$. Let us use $E e^\ld_\ld=0$ applied to the right hand side of (\ref{hwv}) by means of comultiplication $\Delta$. Coefficient near tensor monomial $e^{\ld_1}_{\ld_1-k+1}\ot e^{\ld_2}_{-\ld_1+\ld+k}$ gives the following recursion relation
$$
q^{2(\ld_2-J+k)}\q{k}\qt{\ld_1}{\ld_2}{\ld}{\ld_1-k}{\ld_2-J+k}{\ld}+\q{J-k+1}\qt{\ld_1}{\ld_2}{\ld}{\ld_1-k+1}{\ld_2-J+k-1}{\ld}=0.
$$
Its solution is
$$
\qt{\ld_1}{\ld_2}{\ld}{\ld_1-k}{\ld_2-J+k}{\ld}=\frac{(-1)^k q^{-k(2(\ld_2-J)+k+1)}\p{k-1}{J-j}}{\q{k}!}\qt{\ld_1}{\ld_2}{\ld}{\ld_1}{\ld_2-J}{\ld}.
$$

Now we need to compute
$$
C_{\ld_1,\ld_2,\ld}:=\qt{\ld_1}{\ld_2}{\ld}{\ld_1}{\ld_2-J}{\ld},
$$
which is given by the condition $(e^\ld_\ld,e^\ld_\ld)=1$. We extend the Shapovalov form to tensor products of weight vectors as
$$
(x_1\ot x_2,y_1\ot y_2)=(x_1,y_1)(x_2,y_2)
$$
for any weight vectors $x_i,y_i$, $i=1,2$. Extension defined this way is not functorial, but since the Shapovalov form is not a morphism in the $q$-analog of the category $\mathcal{O}$ anyway, we can proceed with such a definition without bearing any consequences. It is convenient for further computations and influences only the normalization of $q3j$-symbols. We get
\begin{eqnarray}\label{comporth}
&\sum\limits_{k=0}^J\sum\limits_{k^\prime=0}^J\qt{\ld_1}{\ld_2}{\ld}{\ld_1-k}{\ld_2-J+k}{\ld}\qt{\ld_1}{\ld_2}{\ld}{\ld_1-k^\prime}{\ld_2-J+k^\prime}{\ld}(e^{\ld_1}_{\ld_1-k},e^{\ld_1}_{\ld_1-k^\prime})\times\\
&\times(e^{\ld_2}_{\ld_2-J+k},e^{\ld_2}_{\ld_2-J+k^\prime})=\sum\limits_{k=0}^J\left(\qt{\ld_1}{\ld_2}{\ld}{\ld_1-k}{\ld_2-J+k}{\ld}\right)^2\alpha_{\ld_1,k}\alpha_{\ld_2,J-k}=\nonumber\\
&=(C_{\ld_1,\ld_2,\ld})^2\sum\limits_{k=0}^J q^{-2k(2(\ld_2-J)+k+1)}\frac{\p{k-1}{J-j}\p{k-1}{J-j}}{\q{k}!\q{k}!}\frac{(-1)^k q^{k(k-2\ld_1+1)}\q{k}!}{\p{k-1}{2\ld_1-j}}\times\nonumber\\
&\times\frac{(-1)^{J-k}q^{(J-k)(J-k-2\ld_2+1)}\q{J-k}!}{\p{J-k-1}{2\ld_2-j}}=(C_{\ld_1,\ld_2,\ld})^2 q^{\ld(\ld-1)-\ld_1(\ld_1-1)-\ld_2(\ld_2-1)+2\ld_1(\ld_1-\ld)}\times\nonumber\\
&\times(-1)^J \q{J}!\sum\limits_{k=0}^J q^{-2k(\ld+1)}\frac{\p{k-1}{J-j}}{\q{k}!\p{k-1}{2\ld_1-j}\p{J-k-1}{2\ld_2-j}}\boxed{=}.\nonumber
\end{eqnarray}
Using identities (\ref{qid1}) and (\ref{qid2}), the summation in the expression above becomes
\begin{eqnarray*}
(q-q^{-1})^J\sum\limits_{k=0}^J(-1)^k q^{-2k(\ld+1)+\frac{k}{2}(2J+1-k)+\frac{k}{2}(k+1)-\frac{k}{2}(4\ld_1+1-k)-\frac{J-k}{2}(4\ld_2+1-J+k)}\times\\
\times \frac{(q^{-2J};q^2)_k}{(q^2;q^2)_k (q^{-4\ld_1};q^2)_k (q^{-4\ld_2};q^2)_{J-k}}.
\end{eqnarray*}
Then using identity (\ref{id2}) and simplifying further we get
\begin{eqnarray*}
&\boxed{=}(C_{\ld_1,\ld_2,\ld})^2 q^{\ld(\ld-1)-\ld_1(\ld_1-1)-\ld_2(\ld_2-1)+2\ld_1(\ld_1-\ld)+\frac{J}{2}(J-1)-2J\ld_2}\frac{(-1)^J\q{J}! (q-q^{-1})^J}{(q^{-4\ld_2};q^2)_J}\times\\
&\times \sum\limits_{k=0}^\infty q^{-2k(2\ld+1)}\frac{(q^{-2J};q^2)_k (q^{4\ld_2+2(1-J)};q^2)_k}{(q^2;q^2)_k (q^{-4\ld_1};q^2)_k}=(C_{\ld_1,\ld_2,\ld})^2 \frac{(-1)^J\q{J}! (q-q^{-1})^J}{(q^{-4\ld_2};q^2)_J}\times\\
&\times q^{\ld(\ld-1)-\ld_1(\ld_1-1)-\ld_2(\ld_2-1)+2\ld_1(\ld_1-\ld)+\frac{J}{2}(J-1)-2J\ld_2} {}_2 \phi_1\left[\begin{matrix}
q^{-2J}, q^{4\ld_2+2-2J}\\
q^{-4\ld_1}
\end{matrix};q^2,q^{-4\ld-2}\right]=
\end{eqnarray*}
where we took advantage of (\ref{id1}) to change the summation limit to infinity. We use identity (\ref{hid1}) to get
\begin{eqnarray*}
=(C_{\ld_1,\ld_2,\ld})^2 q^{\ld(\ld-1)-\ld_1(\ld_1-1)-\ld_2(\ld_2-1)+2\ld_1(\ld_1-\ld)+\frac{J}{2}(J-1)-2J\ld_2} \times\\
\times (-1)^J\q{J}! (q-q^{-1})^J \frac{(q^{-4\ld_1+2J};q^2)_\infty (q^{-4\ld_1-4\ld_2-2+2J};q^2)_\infty}{(q^{-4\ld_1};q^2)_\infty (q^{-4\ld-2};q^2)_\infty (q^{-4\ld_2};q^2)_J}=
\end{eqnarray*}
and repeatedly making use of (\ref{id1}) to obtain
\begin{eqnarray*}
=(C_{\ld_1,\ld_2,\ld})^2 q^{\ld(\ld-1)-\ld_1(\ld_1-1)-\ld_2(\ld_2-1)+2\ld_1(\ld_1-\ld)+\frac{J}{2}(J-1)-2J\ld_2} \times\\
\times (-1)^J\q{J}! \frac{(q^{-4\ld-2-2J};q^2)_J }{(q^{-4\ld_1};q^2)_J (q^{-4\ld_2};q^2)_J}=
\end{eqnarray*}
then we use (\ref{qid2}) and simplify further to finally arrive at
\begin{eqnarray*}
=(C_{\ld_1,\ld_2,\ld})^2(-1)^J q^{2((\ld-\ld_1)^2-\ld_2^2)}\frac{\q{J}!\p{J-1}{2(\ld+1)+j}}{\p{J-1}{2\ld_1-j}\p{J-1}{2\ld_2-j}}=1.
\end{eqnarray*}
This results in
$$
\qt{\ld_1}{\ld_2}{\ld}{\ld_1}{\ld_2-J}{\ld}=(-1)^{\frac{J}{2}}q^{\ld_2^2-(\ld-\ld_1)^2}\left(\frac{1}{\q{J}!}\prod\limits_{j=0}^{J-1}\frac{\q{2\ld_1-j}\q{2\ld_2-j}}{\q{2\ld+2+j}}\right)^{\frac{1}{2}}.
$$
For $\ld_1+\ld_2-\ld\in\mbb{Z}$ denote
\begin{equation*}
    \dd{\ld_1,\ld_2,\ld}=\left(\q{\ld_1+\ld_2-\ld}!\prod\limits_{j=0}^{\ld_1+\ld_2-\ld-1}\frac{\q{2\ld_1-j}\q{2\ld_2-j}}{\q{2\ld+2+j}}\right)^{\frac{1}{2}}.
\end{equation*}
So, we obtain
\begin{equation*}
\qt{\ld_1}{\ld_2}{\ld}{\ld_1-k}{\ld_2-J+k}{\ld}=\frac{(-1)^{\frac{J}{2}+k}q^{-k(k+2(\ld-\ld_2)+1)+\ld_2^2-(J-\ld_2)^2}}{\q{k}! \q{J-k}!}\dd{\ld_1,\ld_2,\ld}
\end{equation*}
and, hence,
\begin{eqnarray}\label{qthwv}
\qt{\ld_1}{\ld_2}{\ld}{a_1}{a_2}{\ld}=\frac{(-1)^{\frac{J}{2}+\ld_1-a_1}q^{(\ld_1-a_1)((J-\ld_2)-a_2-1)+\ld_2^2-(J-\ld_2)^2}}{\q{\ld_1-a_1}! \q{\ld_2-a_2}!}\dd{\ld_1,\ld_2,\ld}
\end{eqnarray}

\subsubsection{Computing $q3j$-symbols for weight vectors of general case}\label{finalformulaq3j}
Let us use $Fe^\ld_{a+1}=\q{\ld+a+1} e^\ld_a$ applied to (\ref{emb}). Coefficient near tensor monomial $e^{\ld_1}_{a_1}\ot e^{\ld_2}_{a_2}$ gives the following recursion relation
\begin{eqnarray}\label{rec}
   &\q{\ld_1+a_1+1}\qt{\ld_1}{\ld_2}{\ld}{a_1+1}{a_2}{a+1}+q^{-2a_1}\q{\ld_2+a_2+1}\qt{\ld_1}{\ld_2}{\ld}{a_1}{a_2+1}{a+1}=\\
   &=\q{\ld+a+1}\qt{\ld_1}{\ld_2}{\ld}{a_1}{a_2}{a}.\nonumber
\end{eqnarray}
Solution to such recursion is
\begin{eqnarray}\label{ref1}
    &\qt{\ld_1}{\ld_2}{\ld}{a_1}{a_2}{a}=\sum\limits_{k=0}^{\ld-a}q^{-2(\ld-a-k)(a_1+k)+k(\ld-a-k)}\frac{\prod\limits_{j=1}^{k}\q{\ld_1+a_1+j}\prod\limits_{j=1}^{\ld-a-k}\q{\ld_2+a_2+j}}{\prod\limits_{j=1}^{\ld-a}\q{\ld+a+j}}\times\\
    &\times \frac{\q{\ld-a}!}{\q{k}! \q{\ld-a-k}!}\qt{\ld_1}{\ld_2}{\ld}{a_1+k}{a_2+\ld-a-k}{\ld}.\nonumber
\end{eqnarray}
Note that only terms with $(\ld-a)-(\ld_2-a_2)\leq k\leq \ld_1-a_1$ give impact, which is evident from the recursion (\ref{rec}). For example, if we put $a_2=\ld_2$, then only one term $k=\ld-a$ will be present in the summation in (\ref{ref1}). Denote $z=k$, $w=\ld-a-k$ and use (\ref{qthwv}) to arrive at
\begin{eqnarray*}
    \qt{\ld_1}{\ld_2}{\ld}{a_1}{a_2}{a}=(-1)^{\frac{J}{2}+\ld_1-a_1}q^{-(\ld_1-a_1)+\ld_2^2-\ld^2+a^2-a_2^2}\frac{\q{\ld-a}!}{\prod\limits_{j=1}^{\ld-a}\q{\ld+a+j}}\times\\
    \times \dd{\ld_1,\ld_2,\ld}\sum\limits_{z,w}(-1)^zq^{z(\ld+a+1)}\frac{\prod\limits_{j=1}^{z}\q{\ld_1+a_1+j}\prod\limits_{j=1}^{w}\q{\ld_2+a_2+j}}{\q{z}!\q{w}!\q{\ld_1-a_1-z}! \q{\ld_2-a_2-w}!},
\end{eqnarray*}
where $(z,w)\in\adm(\ld-a)$. Again, note that if in any term in at least one of the $q$-factorials the argument is negative, then such term is absent. For example, if $a_2=\ld_2$, the summation in (\ref{rf2}) will have only one term $z=\ld-a$, $w=0$. Also note that for products of $q$-numbers in the formula above we have $\p{-1}{\gamma\pm j}:=1$, $\gamma\in\mbb{C}$ due to initial conditions for recursion relations, which recovers formula (\ref{qthwv}).

\section{Properties of quantum $3j$-symbols}\label{qtprop}
\subsection{Formulae for $q3j$-symbols}
\begin{theorem}\label{theoremq3j}
    $q3j$-symbols in (\ref{emb}) are given by formulae
    \begin{itemize}
        \item (Racah-Fock formula)
    \begin{eqnarray}\label{rf2}
    \qt{\ld_1}{\ld_2}{\ld}{a_1}{a_2}{a}^\psi=(-1)^{\frac{J}{2}+\ld_1-a_1}q^{(\ld_2^2-a_2^2)-(\ld^2-a^2)-(\ld_1-a_1)}\frac{\q{\ld-a}!}{\prod\limits_{j=1}^{\ld-a}\q{\ld+a+j}}\times\\
    \times \dd{\ld_1,\ld_2,\ld}\sum\limits_{z,w}(-1)^zq^{z(\ld+a+1)}\frac{\prod\limits_{j=1}^{z}\q{\ld_1+a_1+j}\prod\limits_{j=1}^{w}\q{\ld_2+a_2+j}}{\q{z}!\q{w}!\q{\ld_1-a_1-z}! \q{\ld_2-a_2-w}!},\nonumber
    \end{eqnarray}
    where $(z,w)\in\adm(\ld-a)$.
        \item (Van der Waerden formula)
    \begin{eqnarray}\label{vdw2}
        &\qt{\ld_1}{\ld_2}{\ld}{a_1}{a_2}{a}^\psi=(-1)^{\frac{J}{2}}q^{\ld_1^2+\ld_2^2-\ld^2+\ld_1\ld_2+a_1a_2+\ld_1a_2-\ld_2a_1}\frac{\q{\ld-a}!}{\prod\limits_{j=1}^{\ld-a}\q{\ld+a+j}}\times\\
        &\times \dd{\ld_1,\ld_2,\ld}\sum\limits_{z,w}(-1)^z q^{-z(\ld_1+\ld_2+\ld+1)}\frac{\q{\ld_2-\ld_1+\ld|\ld_2+a_2-z}\q{\ld_1-\ld_2+\ld|\ld_1+a_1-w}}{\q{z}!\q{w}!\q{\ld_1-a_1-z}! \q{\ld_2-a_2-w}!},\nonumber
    \end{eqnarray}
    where $(z,w)\in\adm(J)$.
    \end{itemize}
\end{theorem}
Function $\q{a|b}$ is defined in (\ref{qpf}). There is also an analog of Majumdar formula for $q3j$-symbols for Verma modules, which is given by a change of variables $z\to \ld-\ld_2-a_1+z$. In what follows we won't need it.
\begin{proof}
    ({\it Racah-Fock formula})\\
    Expression (\ref{rf2}) has been obtained as a result of Section \ref{deriveformulaq3j}.

    ({\it Van der Waerden formula})\\
    Consider the equality between the Van der Waerden formula and the Racah-Fock formula in the finite-dimensional case \cite{K}. Multiply both sides of the equality by $\sqrt{\frac{\q{-\ld_1+\ld_2+\ld}!\q{\ld_1-\ld_2+\ld}!}{\q{\ld_1+a_1}!\q{\ld_2+a_2}!}}$. This results in the following identity
    \begin{eqnarray}\label{ididid}
        &(-1)^{\ld_1-a_1}\sum\limits_{z,w}(-1)^zq^{z(\ld+a+1)}\frac{\prod\limits_{j=1}^{z}\q{\ld_1+a_1+j}\prod\limits_{j=1}^{w}\q{\ld_2+a_2+j}}{\q{z}!\q{w}!\q{\ld_1-a_1-z}! \q{\ld_2-a_2-w}!}=\\
        &=q^{(\ld_1-a_1)(\ld_2+a_2+1)+\ld_1^2-a_1^2}\sum\limits_{u,v}(-1)^u q^{-u(\ld_1+\ld_2+\ld+1)}\frac{\q{\ld_2-\ld_1+\ld|\ld_2+a_2-u}\q{\ld_1-\ld_2+\ld|\ld_1+a_1-v}}{\q{u}!\q{v}!\q{\ld_1-a_1-u}! \q{\ld_2-a_2-v}!},\nonumber
    \end{eqnarray}
    where $(z,w)\in\adm(\ld-a)$ and $(u,v)\in\adm(J)$. Both sides of the equality we initially used are $q$-orthogonal polynomials, so the identity (\ref{ididid}) admits analytic continuation to $\ld_1,\ld_2,\ld\in\mbb{C}$. The resulting formula is well-defined since all the $q$-factorials have arguments in $\mbb{Z}_{\geq 0}$. Plugging this identity into formula (\ref{rf2}) gives the desired result.
\end{proof}
From (\ref{prj}) and computation in (\ref{comporth}) we deduce that
\begin{equation}\label{projdeduce}
    \qt{\ld_1}{\ld_2}{\ld}{a_1}{a_2}{a}^\pi=\frac{\alpha_{\ld_1,\ld_1-a_1}\alpha_{\ld_2,\ld_2-a_2}}{\alpha_{\ld,\ld-a}}\qt{\ld_1}{\ld_2}{\ld}{a_1}{a_2}{a}^\psi.
\end{equation}

\subsection{String diagrams for $q3j$-symbols}\label{strings}
\begin{figure}[h]
\begin{center}
	{\includegraphics[width=420pt]{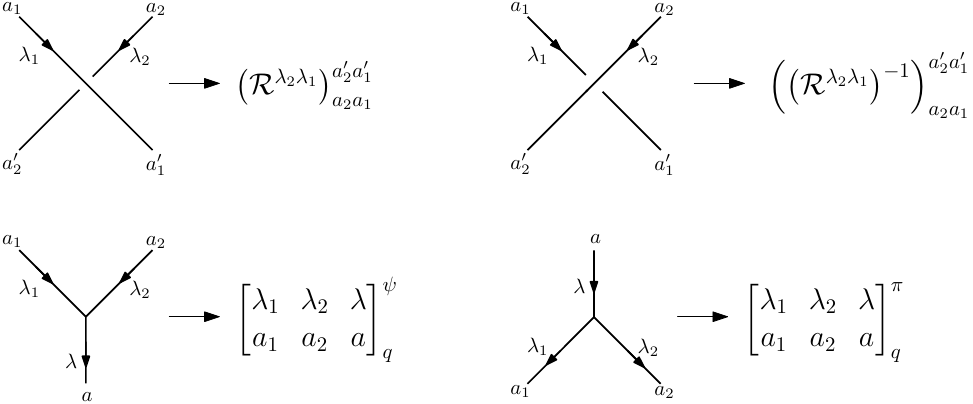}}
    \caption{Assignment of string diagrams for $\mathcal{R}$-matrix and $q3j$-symbols.}
\label{assign}
\end{center}
\end{figure}

Following \cite{KR}, we will examine properties of $q3j$-symbols using string diagrams. Strings are colored by highest weights $\ld_j$ of Verma modules $M_{\ld_j}$, their ends are marked by weights $a_j$ of vectors $e^{\ld_j}_{a_j}$ in weight spaces of representations coloring the strings. Crossings correspond to the $\mathcal{R}$-matrix and its inverse, each trivalent vertex with two ends joining in one corresponds to an embedding $\psi$, each trivalent vertex with one end splitting in two corresponds to a projection $\pi$. This assignment is depicted in Figure \ref{assign}.
\par
\begin{figure}[h]
\begin{center}
	{\includegraphics[width=440pt]{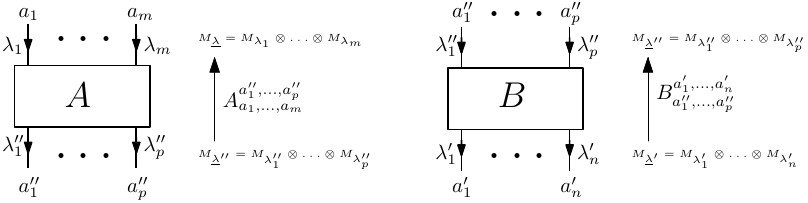}}
    \caption{String diagrams for morphisms $A:M_{\underline{\ld}^{\prime\prime}}\to M_{\underline{\ld}}$ and $B:M_{\underline{\ld}^{\prime}}\to M_{\underline{\ld}^{\prime\prime}}$, where $A_{a_1,\ldots,a_m}^{a^{\prime\prime}_1,\ldots,a^{\prime\prime}_p}$ and $B^{a_1^\prime,\ldots,a_n^\prime}_{a^{\prime\prime}_1,\ldots,a^{\prime\prime}_p}$ are corresponding matrix elements.}
    \label{assign2}
\end{center}
\end{figure}

Two morphisms of a general form $A:M_{\underline{\ld}^{\prime\prime}}\to M_{\underline{\ld}}$, $B:M_{\underline{\ld}^{\prime}}\to M_{\underline{\ld}^{\prime\prime}}$ can be represented graphically by string diagrams as in Figure \ref{assign2}.
\par
\begin{figure}[h]
\begin{center}
	{\includegraphics[width=380pt]{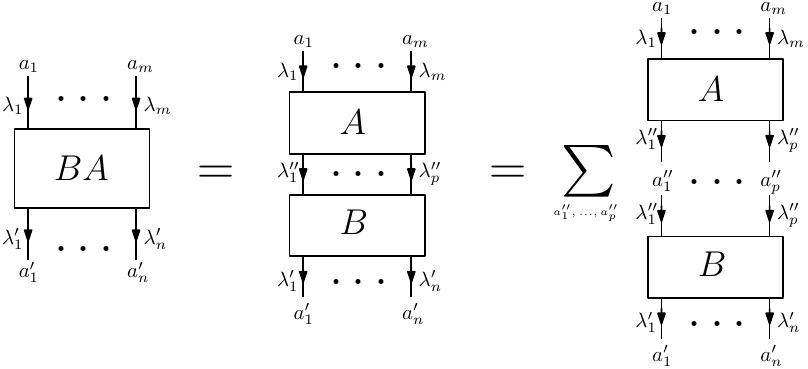}}
    \caption{String diagram for composition of morphisms $A$ and $B$.}
    \label{assign3}
\end{center}
\end{figure}

The composition of these morphisms is represented by joining the ends of the strings of the corresponding graphs and summing over all weights at the ends of the joining strings, as shown in Figure \ref{assign3}. This procedure amounts to the summation of products of matrix elements $B^{a_1^\prime,\ldots,a_n^\prime}_{a^{\prime\prime}_1,\ldots,a^{\prime\prime}_p}A_{a_1,\ldots,a_m}^{a^{\prime\prime}_1,\ldots,a^{\prime\prime}_p}$ over all values of weights $a^{\prime\prime}_j$, where $j=1,\ldots,p$.
\par
\begin{figure}[h]
\begin{center}
	{\includegraphics[width=200pt]{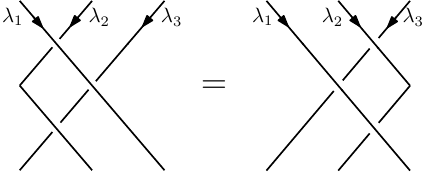}}
    \caption{String diagram representation for the Yang-Baxter relation (\ref{YB}).}
    \label{q3jYB}
\end{center}
\end{figure}

For future convenience we will omit writing marks at the ends of the strings, implying the composition rule above. For example, Yang-Baxter relation (\ref{YB}) is represented graphically in Figure \ref{q3jYB}.
\par
Note that we are not considering graphical representation of the Shapovalov form since it is not functorial and does not respect the topology of string diagrams. We also restrict ourselves to considering only strings, as opposed to ribbon graphs, since definition of the ribbon structure should also involve Verma modules of the lowest weight. This problem will be discussed in Section \ref{chevv}.

\subsection{Orthogonality and completeness relations}
\begin{theorem}\label{orthrelth}
    The projection $\pi_\ld:M_{\ld_1}\otimes M_{\ld_2}\to M_\ld$ is $\Uq$-invariant and
    \begin{equation*}
        \pi_\ld \circ \psi_\ld = \mathrm{Id}_{M_\ld},\quad \sum\limits_{\ld} \psi_\ld \circ \pi_\ld = \mathrm{Id}_{M_{\ld_1}\otimes M_{\ld_2}},
    \end{equation*}
    or, explicitly,
    \begin{equation}\label{projth1}
    \sum\limits_{a_1,a_2}\qt{\ld_1}{\ld_2}{\ld}{a_1}{a_2}{a}^\pi\qt{\ld_1}{\ld_2}{\ld^\prime}{a_1}{a_2}{a^\prime}^\psi=\delta_{\ld,\ld^\prime} \delta_{a,a^\prime},
    \end{equation}
\begin{figure}[h]
\begin{center}
	{\includegraphics[width=150pt]{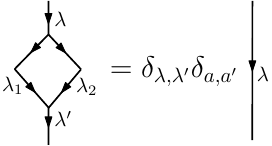}}
    \caption{String diagram representations for the orthagonality relation (\ref{projth1}).}
\end{center}
\end{figure}\par
    \begin{equation}\label{projth2}
    \sum\limits_{\ld,a}\qt{\ld_1}{\ld_2}{\ld}{a_1}{a_2}{a}^\pi\qt{\ld_1}{\ld_2}{\ld}{a_1^\prime}{a_2^\prime}{a}^\psi=\delta_{a_1,a_1^\prime} \delta_{a_2,a_2^\prime}.
    \end{equation}
\begin{figure}[h]
\begin{center}
	{\includegraphics[width=230pt]{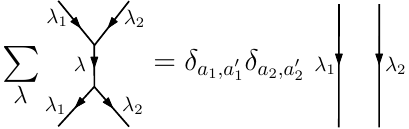}}
    \caption{String diagram representations for the completeness relation (\ref{projth2}).}
\end{center}
\end{figure}\par
\end{theorem}
\begin{proof}
    It is straightforward to check $\Uq$-invariance, it follows from recursion relations (\ref{rec}) and (\ref{recmor}). Consider condition $(e^\ld_a,e^{\ld^\prime}_{a^\prime})=\delta_{a,a^\prime}\delta_{\ld,\ld^\prime}\alpha_{\ld,\ld-a}$. Computation similar to (\ref{comporth}) gives the result. The same idea is used for proving (\ref{projth2}).
\end{proof}
Note that in (\ref{projth2}) once $a_1$, $a_2$ and $a_1^\prime$, $a_2^\prime$ are fixed, so are $a$ and $a^\prime$, and summation over $\ld$ becomes finite since $\ld-a\geq 0$ and $\ld_1+\ld_2-\ld\geq 0$.

\subsection{Identities for $q3j$-symbols and the universal $\mathcal{R}$-matrix}
\begin{theorem}\label{rmatth}
The following identities hold
\begin{eqnarray}\label{rmatid1}
    &\sum\limits_{a_1^\prime,a_2^\prime}\left(\mathcal{R}^{\ld_1,\ld_2}\right)^{a_1^\prime,a_2^\prime}_{a_1,a_2}\qt{\ld_1}{\ld_2}{\ld}{a_1^\prime}{a_2^\prime}{a}^\psi=(-1)^Jq^{C-C_1-C_2}\qt{\ld_2}{\ld_1}{\ld}{a_2}{a_1}{a}^\psi,
\end{eqnarray}
where
$$
C_j:=\ld_j(\ld_j+1),\quad C:=\ld(\ld+1),
$$
\begin{figure}[h]
\begin{center}
	{\includegraphics[width=265pt]{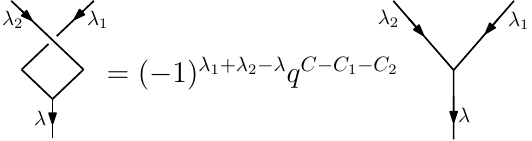}}
    \caption{String diagram representations for the relation (\ref{rmatid1}).}
\end{center}
\end{figure}\par
and
\begin{eqnarray}\label{rmatid2}
    &\sum\limits_{a^\prime}\left(\mathcal{R}^{\ld,\ld_3}\right)^{a,a_3^\prime}_{a^\prime,a_3}\qt{\ld_1}{\ld_2}{\ld}{a_1}{a_2}{a^\prime}^\psi=\sum\limits_{a_1^\prime,a_2^\prime,a_3^{\prime\prime}}\left(\mathcal{R}^{\ld_1,\ld_3}\right)^{a_1^\prime,a_3^{\prime\prime}}_{a_1,a_3}\left(\mathcal{R}^{\ld_2,\ld_3}\right)^{a_2^\prime,a_3^\prime}_{a_2,a_3^{\prime\prime}}\qt{\ld_1}{\ld_2}{\ld}{a_1^\prime}{a_2^\prime}{a}^\psi.
\end{eqnarray}
\begin{figure}[h]
\begin{center}
	{\includegraphics[width=205pt]{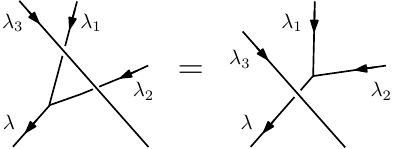}}
    \caption{String diagram representations for the relation (\ref{rmatid1}).}
\end{center}
\end{figure}\par
\end{theorem}
\begin{proof}
The proof of these identities is the same as for finite-dimensional representations. For details see \cite{R1}, Theorem 1.5 and Theorem 1.6 for (\ref{rmatid1}) and (\ref{rmatid2}) respectively.
\end{proof}
Note that in (\ref{rmatid1}) summations in both sides come from contractions. But we have a condition $a^\prime=a_1+a_2$, which comes from $q3j$-symbol in the left hand side. This combined with formula (\ref{rmatrel}) gives $a_3^\prime=a_3$. From these conditions we also get $a_3^{\prime\prime}=a_3$, $a_1^\prime=a_1$, $a_2^\prime=a_2$. As a result, formula (\ref{rmatid2}) is reduced to
\begin{eqnarray*}
    \left(\mathcal{R}^{\ld,\ld_3}\right)^{a,a_3}_{a,a_3}\qt{\ld_1}{\ld_2}{\ld}{a_1}{a_2}{a}^\psi=\left(\mathcal{R}^{\ld_1,\ld_3}\right)^{a_1,a_3}_{a_1,a_3}\left(\mathcal{R}^{\ld_2,\ld_3}\right)^{a_2,a_3}_{a_2,a_3}\qt{\ld_1}{\ld_2}{\ld}{a_1}{a_2}{a}^\psi.
\end{eqnarray*}
Identities (\ref{rmatid1}) and (\ref{rmatid2}) are also true for projections.

\subsection{Symmetries of $q3j$-symbols}\label{qtsectsym}
\begin{theorem}
The following identities hold
\begin{equation*}
    \qt{\ld_1}{\ld_2}{\ld}{a_1}{a_2}{a}^\psi=(-1)^{-J}q^{(\ld_1-a_1)(\ld_1+a_1-1)+(\ld_2-a_2)(\ld_2+a_2-1)-(\ld-a)(\ld+a-1)}\qtm{\ld_2}{\ld_1}{\ld}{a_2}{a_1}{a}^\psi,
\end{equation*}
\begin{equation*}
    \qt{\ld_1}{\ld_2}{\ld}{a_1}{a_2}{a}^\pi=(-1)^{-J}q^{-(\ld_1-a_1)(\ld_1+a_1-1)-(\ld_2-a_2)(\ld_2+a_2-1)+(\ld-a)(\ld+a-1)}\qtm{\ld_2}{\ld_1}{\ld}{a_2}{a_1}{a}^\pi.
\end{equation*}
\end{theorem}
\begin{proof}
    Follows from a straightforward computation using (\ref{rf2}). Identity for projections follows from the identity for embeddings and Theorem \ref{orthrelth}.
\end{proof}
The results obtained here extend naturally to tensor products involving finite-dimensional irreducible representations \cite{BG,BGG}, provided the highest weights of the components are chosen generically so that all direct summands in the decomposition remain irreducible. In this setting, the previously derived formulae continue to apply, subject to appropriate restrictions on the admissible values of $J$ and on the weights $a_j$ corresponding to the finite-dimensional components. Note that $q3j$-symbols for Verma modules admit fewer expressions than their finite-dimensional counterparts, which possess additional symmetries from self-duality.
\section{Quantum $6j$-symbols for Verma modules}\label{sectionq6j}
We will proceed with derivation of $q6j$-symbols, also known as quantum Racah-Wigner coefficients. The results of this section are summarized in Section \ref{sectionpropq6j}.
\subsection{Definition of the $q6j$-symbols}
Consider $M_{\ld_1}\ot M_{\ld_2}\ot M_{\ld_3}$, which can be decomposed in two ways. As $(M_{\ld_1}\ot M_{\ld_2})\ot M_{\ld_3}$:
\begin{equation}\label{dec1}
    e^{\ld_{12},\ld}_a=\sum\limits_{a_1,a_2,a_3}\qt{\ld_{12}}{\ld_3}{\ld}{a_{12}}{a_3}{a}^\psi\qt{\ld_1}{\ld_2}{\ld_{12}}{a_1}{a_2}{a_{12}}^\psi e^{\ld_1}_{a_1}\ot e^{\ld_2}_{a_2}\ot e^{\ld_3}_{a_3},
\end{equation}
where $(a_{12},a_3)\in\adm(\ld_{12},\ld_3,\ld_{12}+\ld_3-a)$, $(a_1,a_2)\in\adm(\ld_1,\ld_2,\ld_1+\ld_2-a_{12})$. And as $M_{\ld_1}\ot (M_{\ld_2}\ot M_{\ld_3})$:
\begin{equation}\label{dec2}
    e^{\ld_{23},\ld}_a=\sum\limits_{a_1,a_2,a_3}\qt{\ld_1}{\ld_{23}}{\ld}{a_1}{a_{23}}{a}^\psi\qt{\ld_2}{\ld_3}{\ld_{23}}{a_2}{a_3}{a_{23}}^\psi e^{\ld_1}_{a_1}\ot e^{\ld_2}_{a_2}\ot e^{\ld_3}_{a_3},
\end{equation}
where $(a_1,a_{23})\in\adm(\ld_1,\ld_{23},\ld_1+\ld_{23}-a)$, $(a_2,a_3)\in\adm(\ld_2,\ld_3,\ld_2+\ld_3-a_{23})$. Also, 
$$
\ld_1+\ld_2-\ld_{12},\ld_2+\ld_3-\ld_{23}, \ld_1+\ld_2+\ld_3-\ld\in\mbb{Z}_{\geq 0},
$$ 
and 
$$
\ld_1+\ld_2-\ld_{12}\leq \ld_1+\ld_2+\ld_3-\ld,\quad \ld_2+\ld_3-\ld_{23}\leq \ld_1+\ld_2+\ld_3-\ld.
$$
\begin{definition}
    Matrix elements of the linear transformation relating decompositions (\ref{dec1}) and (\ref{dec2})
    \begin{equation}\label{qsdef}
        e^{\ld_{12},\ld}_a=\sum\limits_{\ld_{23}}\qs{\ld_1}{\ld_2}{\ld_{12}}{\ld_3}{\ld}{\ld_{23}} e^{\ld_{23},\ld}_a
    \end{equation}
    are called $q6j$-symbols.
\begin{figure}[h]
\begin{center}
	{\includegraphics[width=350pt]{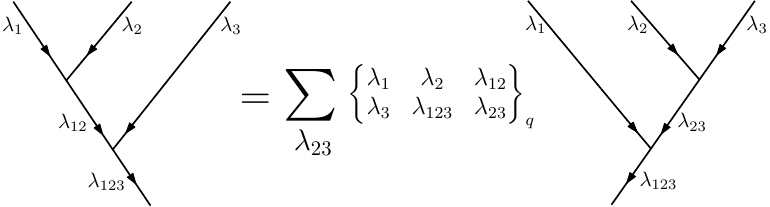}}
    \caption{String diagram representation for the definition of $q6j$-symbols (\ref{qsdef}), which relates identities (\ref{dec1}) and (\ref{dec2}).}
\end{center}
\end{figure}\par
\end{definition}
Note that summation over $\ld_{23}$ is finite. For fixed $\ld_1$, $\ld_2$, $\ld_3$, $\ld$ from decomposition (\ref{decomp}) we have $\ld_1+\ld_{23}-\ld\geq 0$ and $\ld_2+\ld_3-\ld_{23}\geq 0$.

\subsection{Computing $q6j$-symbols}\label{deriveformulaq6j}
In the following subsections, we apply the Racah method \cite{BL} to derive an explicit formula for the $q6j$-symbols. For clarity, we break the procedure into steps, with each subsection corresponding to a single step.
\subsubsection{Step 0}
In this step we provide a starting point for the Racah method. Consider (\ref{qsdef}) and substitute (\ref{dec1}) and (\ref{dec2})
\begin{eqnarray}\label{step01}
    \sum\limits_{a_1,a_2,a_3}\qt{\ld_{12}}{\ld_3}{\ld}{a_{12}}{a_3}{a}^\psi\qt{\ld_1}{\ld_2}{\ld_{12}}{a_1}{a_2}{a_{12}}^\psi=\sum\limits_{\ld_{23}}\qs{\ld_1}{\ld_2}{\ld_{12}}{\ld_3}{\ld}{\ld_{23}}\times\\
    \times\sum\limits^\prime_{a_1,a_2,a_3}\qt{\ld_1}{\ld_{23}}{\ld}{a_1}{a_{23}}{a}^\psi\qt{\ld_2}{\ld_3}{\ld_{23}}{a_2}{a_3}{a_{23}}^\psi\nonumber,
\end{eqnarray}
where we have omitted writing $e^{\ld_1}_{a_1}\ot e^{\ld_2}_{a_2}\ot e^{\ld_3}_{a_3}$. Summations over $a_1,a_2,a_3$ on both sides are different, in accordance to (\ref{dec1}) and (\ref{dec2}). We apply projection
\begin{equation*}
    \pi_{\ld_{23}}(e^{\ld_2}_{a_2}\ot e^{\ld_3}_{a_3})=\qt{\ld_2}{\ld_3}{\ld_{23}}{a_2}{a_3}{a_2+a_3}^\pi e^{\ld_{23}}_{a_2+a_3}
\end{equation*}
to both sides of (\ref{step01}) and use orthogonality relation (\ref{projth1}), which gets rid of the summation over $\ld_{23}$, $a_2$, $a_3$ on the right hand side. Since $a_{23}$ is now also fixed, summation over $a_1$ is gone as well. This allows us to divide both sides of the obtained formula by
\begin{equation*}
\qt{\ld_1}{\ld_{23}}{\ld}{a_1}{a_{23}}{a}^\psi
\end{equation*}
and use (\ref{projdeduce}) to get
\begin{eqnarray}\label{step02}
    \qs{\ld_1}{\ld_2}{\ld_{12}}{\ld_3}{\ld}{\ld_{23}}=\left(\qt{\ld_1}{\ld_{23}}{\ld}{a_1}{a_{23}}{a}^\psi\right)^{-1}\sum\limits_{a_1,a_2,a_3}\frac{\alpha_{\ld_2,\ld_2-a_2}\alpha_{\ld_3,\ld_3-a_3}}{\alpha_{\ld_{23},\ld_{23}-a_{23}}}\times\\
    \times \qt{\ld_2}{\ld_3}{\ld_{23}}{a_2}{a_3}{a_{23}}^\psi\qt{\ld_{12}}{\ld_3}{\ld}{a_{12}}{a_3}{a}^\psi\qt{\ld_1}{\ld_2}{\ld_{12}}{a_1}{a_2}{a_{12}}^\psi\nonumber,
\end{eqnarray}
where $(a_{12},a_3)\in\adm(\ld_{12},\ld_3,\ld_{12}+\ld_3-a)$, $(a_1,a_2)\in\adm(\ld_1,\ld_2,\ld_1+\ld_2-a_{12})$. From now on we omit writing embedding index $\psi$ above $q3j$-symbols. Since the left hand side of (\ref{step02}) is independent of some of the parameters present on the right hand side, we can put
\begin{equation*}
    a\to\ld,\quad a_{23}\to\ld_{23},\quad a_3\to\ld_3-k,
\end{equation*}
which fixes other parameters as
\begin{equation*}
    a_1\to\ld-\ld_{23},\quad a_2\to\ld_{23}-\ld_3+k,\quad a_{12}\to\ld-\ld_3+k,
\end{equation*}
and from the range of $a_3$ we have 
\begin{equation*}
   k=0,1,\dots,\ld_{12}+\ld_3-\ld. 
\end{equation*}
Hence, we obtain
\begin{eqnarray*}
    &\qs{\ld_1}{\ld_2}{\ld_{12}}{\ld_3}{\ld}{\ld_{23}}=\left(\qt{\ld_1}{\ld_{23}}{\ld}{\ld-\ld_{23}}{\ld_{23}}{\ld}\right)^{-1}\sum\limits_{k=0}^{\ld_{12}+\ld_3-\ld}\alpha_{\ld_2,\ld_2+\ld_3-\ld_{23}-k}\alpha_{\ld_3,k}\times\\
    &\times \qt{\ld_2}{\ld_3}{\ld_{23}}{\ld_{23}-\ld_3+k}{\ld_3-k}{\ld_{23}}\qt{\ld_{12}}{\ld_3}{\ld}{\ld-\ld_3+k}{\ld_3-k}{\ld}\qt{\ld_1}{\ld_2}{\ld_{12}}{\ld-\ld_{23}}{\ld_{23}-\ld_3+k}{\ld-\ld_3+k}.
\end{eqnarray*}
Now use (\ref{alphaa}), (\ref{qthwv}) and (\ref{vdw2}) to get
\begin{eqnarray*}
    &\qs{\ld_1}{\ld_2}{\ld_{12}}{\ld_3}{\ld}{\ld_{23}}=(-1)^{\ld_1+\ld_2-\ld_{12}}q^{\ld_1(\ld_1+1)+\ld_3(\ld_3+1)-\ld_{12}(\ld_{12}+1)+\ld_{23}(\ld_1+\ld_2+\ld+1)+\ld_2(\ld_1-\ld)}\times\\
    &\times \dd{\ld_1,\ld_{23},\ld}\dd{\ld_2,\ld_3,\ld_{23}}\dd{\ld_{12},\ld_3,\ld}\dd{\ld_1,\ld_2,\ld_{12}}\prod\limits_{j=0}^{\ld_1+\ld_{23}-\ld-1}\frac{\q{2\ld+2+j}}{\q{2\ld_1-j}\q{2\ld_{23}-j}}\times\\
    &\times \sum\limits_{k=0}^{\ld_3+\ld_{12}-\ld}\frac{q^{-(\ld_3-k)(\ld_1+\ld+\ld_{23}+2)}}{\q{k}!\prod\limits_{j=0}^{k-1}\q{2\ld_3-j}\prod\limits_{j=0}^{\ld_2+\ld_3-\ld_{23}-k-1}\q{2\ld_2-j}\prod\limits_{j=0}^{\ld_3+\ld_{12}-\ld-k-1}\q{2\ld_{12}-j}}\times\\
    &\times \sum\limits_{z=0}^{\ld_1+\ld_2-\ld_{12}}(-1)^z q^{-z(\ld_1+\ld_2+\ld_{12}+1)}\frac{\q{\ld_1-\ld_2+\ld_{12}|\ld+\ld_{12}-\ld_2-\ld_{23}+z}\q{\ld_2-\ld_1+\ld_{12}|\ld_2+\ld_{23}-\ld_3+k-z}}{\q{z}!\q{\ld_1+\ld_2-\ld_{12}-z}!\q{\ld_1+\ld_{23}-\ld-z}! \q{\ld_{12}+\ld_3-\ld_1-\ld_{23}-k+z}!}
\end{eqnarray*}

\subsubsection{Step 1}
Consider summation over $k$ and denote $\ld_3-k=\al$.
\begin{eqnarray*}
    &\sum\limits_{\al=\ld-\ld_{12}}^{\ld_3}\frac{q^{-\al(\ld_1+\ld+\ld_{23}+2)}\q{\ld_2-\ld_1+\ld_{12}|\ld_2+\ld_{23}-\al-z}}{\q{\ld_3-\al}!\q{\ld_{12}-\ld_1-\ld_{23}+\al+z}!\prod\limits_{j=0}^{\ld_3-\al-1}\q{2\ld_3-j}\prod\limits_{j=0}^{\ld_2-\ld_{23}+\al-1}\q{2\ld_2-j}\prod\limits_{j=0}^{\ld_{12}-\ld+\al-1}\q{2\ld_{12}-j}}=\\
    &=\sum\limits_{\al=\ld-\ld_{12}}^{\ld_3}\sum\limits_{t:t-\ld_3\in\mbb{Z}}\frac{\delta_{t,\al}q^{-\al(\ld_1+\ld+\ld_{23}+2)}\q{\ld_2-\ld_1+\ld_{12}|\ld_2+\ld_{23}-\al-z}}{\q{\ld_3-\al}!\q{\ld_{12}-\ld_1-\ld_{23}+\al+z}!\prod\limits_{j=0}^{\ld_3-\al-1}\q{2\ld_3-j}\prod\limits_{j=0}^{\ld_2-\ld_{23}+\al-1}\q{2\ld_2-j}\prod\limits_{j=0}^{\ld_{12}-\ld+\al-1}\q{2\ld_{12}-j}}=\\
    &=\sum\limits_{\al=\ld-\ld_{12}}^{\ld_3}\sum\limits_{t:t-\ld_3\in\mbb{Z}}\frac{\delta_{t,\al}q^{-\al(\ld_1+\ld+\ld_{23}+2)+{\color{red}(t-\al)(\ld_3+\ld_2+\ld_{23}-\ld-\ld_{12}-z-1)}}\q{\ld_2-\ld_1+\ld_{12}|\ld_2+\ld_{23}-{\color{red} t}-z}}{\q{\ld_3-{\color{red} t}}!\q{\ld_{12}-\ld_1-\ld_{23}+\al+z}!\prod\limits_{j=0}^{\ld_3-\al-1}\q{2\ld_3-j}\prod\limits_{j=0}^{\ld_2-\ld_{23}+\al-1}\q{2\ld_2-j}\prod\limits_{j=0}^{\ld_{12}-\ld+{\color{red} t}-1}\q{2\ld_{12}-j}}=,
\end{eqnarray*}
where in {\color{red} red} we have highlighted changes we made under the Kronecker delta $\delta_{t,\al}$. Now we use (\ref{qbid1}) to arrive at
\begin{eqnarray*}
    &=\sum\limits_{\al=\ld-\ld_{12}}^{\ld_3}\sum\limits_{\substack{s,t:\\ t-\ld_3\in\mbb{Z}\\ s-\ld_3\in\mbb{Z}}}\frac{(-1)^{t-s}q^{-s(t-\al)+(s-\al)}q^{-\al(\ld_1+\ld+\ld_{23}+2)+(t-\al)(\ld_3+\ld_2+\ld_{23}-\ld-\ld_{12}-z-1)}\q{\ld_2-\ld_1+\ld_{12}|\ld_2+\ld_{23}-t-z}}{\q{t-s}!\q{s-\al}!\q{\ld_3- t}!\q{\ld_{12}-\ld_1-\ld_{23}+\al+z}!\prod\limits_{j=0}^{\ld_3-\al-1}\q{2\ld_3-j}\prod\limits_{j=0}^{\ld_2-\ld_{23}+\al-1}\q{2\ld_2-j}\prod\limits_{j=0}^{\ld_{12}-\ld+t-1}\q{2\ld_{12}-j}}.
\end{eqnarray*}
Note that summation over $s,t$ is finite, since $t$ can have values at most $\ld-\ld_{12},\dots,\ld_3$ and $s=\al,\dots, t$.

\subsubsection{Step 2}
In this step we will get rid of the summation over $\al$. Consider
\begin{eqnarray*}
    &\sum\limits_{\al}\frac{q^{-\al(\ld_3+\ld_2+\ld_1+2\ld_{23}-\ld_{12}-z-s+2)}}{\q{s-\al}!\q{\ld_{12}-\ld_1-\ld_{23}+\al+z}!\prod\limits_{j=0}^{\ld_3-\al-1}\q{2\ld_3-j}\prod\limits_{j=0}^{\ld_2-\ld_{23}+\al-1}\q{2\ld_2-j}}=q^{-s(\ld_3+\ld_2+\ld_1+2\ld_{23}-\ld_{12}-z-s+2)}\times\\
    &\times\sum\limits_{\al}\frac{q^{(s-\al)(\ld_3+\ld_2+\ld_1+2\ld_{23}-\ld_{12}-z-s+2)}}{\q{s-\al}!\q{\ld_{12}-\ld_1-\ld_{23}+s+z-(s-\al)}!\prod\limits_{j=0}^{\ld_3-s+(s-\al)-1}\q{2\ld_3-j}\prod\limits_{j=0}^{\ld_2-\ld_{23}+s-(s-\al)-1}\q{2\ld_2-j}}=,
\end{eqnarray*}
denote $s-\al=p$ and use Lemma \ref{lemstep2} to get
\begin{eqnarray*}
    =q^{-s(\ld_3+\ld_2+\ld_1+2\ld_{23}-\ld_{12}-z-s+2)}q^{(\ld_2+\ld_{23}-s+1)(\ld_{12}-\ld_1-\ld_{23}+z+s)}\times\\
    \times\frac{\q{\ld_2+\ld_3+\ld_{23}+1|\ld_1+\ld_2+\ld_3+2\ld_{23}-\ld_{12}-z-s+1}}{\q{\ld_{12}-\ld_1-\ld_{23}+z+s}!\prod\limits_{j=0}^{\ld_3-\ld_1-\ld_{23}+\ld_{12}+z-1}\q{2\ld_3-j}\prod\limits_{j=0}^{\ld_2-\ld_{23}+s-1}\q{2\ld_2-j}}.
\end{eqnarray*}

\subsubsection{Step 3}
In this step we will get rid of the summation over $t$. Consider
\begin{eqnarray*}
    &\sum\limits_{t}(-1)^t q^{t(\ld_3+\ld_2+\ld_{23}-\ld-\ld_{12}-z-s-1)}\frac{\q{\ld_2-\ld_1+\ld_{12}|\ld_2+\ld_{23}-t-z}}{\q{\ld_3-t}!\q{t-s}!\prod\limits_{j=0}^{\ld_{12}-\ld+t-1}\q{2\ld_{12}-j}}=(-1)^s q^{s(\ld_3+\ld_2+\ld_{23}-\ld-\ld_{12}-z-s-1)}\times\\
    &\times\sum\limits_{t}(-1)^{t-s} q^{(t-s)(\ld_3+\ld_2+\ld_{23}-\ld-\ld_{12}-z-s-1)}\frac{\q{\ld_2-\ld_1+\ld_{12}|\ld_2+\ld_{23}-s-z-(t-s)}}{\q{\ld_3-s-(t-s)}!\q{t-s}!\prod\limits_{j=0}^{\ld_{12}-\ld+s+(t-s)-1}\q{2\ld_{12}-j}}=,
\end{eqnarray*}
denote $t-s=p$ and use Lemma \ref{lemstep3} to get
\begin{eqnarray*}
    =(-1)^s q^{s(\ld_3+\ld_2+\ld_{23}-\ld-\ld_{12}-z-s-1)}q^{(\ld_3-s)(\ld_2+\ld_{23}-z-s)}\times\\
    \times \frac{\q{\ld_2-\ld_1+\ld_{12}|\ld_2+\ld_{23}-s-z}\prod\limits_{j=0}^{\ld_3-s-1}\q{\ld+\ld_{12}-\ld_2-\ld_{23}+z-j}}{\q{\ld_3-s}!\prod\limits_{j=0}^{\ld_{12}+\ld_3-\ld-1}\q{2\ld_{12}-j}}.
\end{eqnarray*}

\subsubsection{Intermediate result after Steps 1-3}
From (\ref{qpfid2}) it follows that
\begin{eqnarray*}
    &\q{\ld_1-\ld_2+\ld_{12}|\ld+\ld_{12}-\ld_2-\ld_{23}+z}\prod\limits_{j=0}^{\ld_3-s-1}\q{\ld+\ld_{12}-\ld_2-\ld_{23}+z-j}=\\
    &=\q{\ld_1-\ld_2+\ld_{12}|\ld+\ld_{12}-\ld_2-\ld_3-\ld_{23}+z+s}.
\end{eqnarray*}
Hence, we obtain
\begin{eqnarray*}
    \qs{\ld_1}{\ld_2}{\ld_{12}}{\ld_3}{\ld}{\ld_{23}}=(-1)^{\ld_1+\ld_2-\ld_{12}}q^{\ld_{12}(\ld_2+\ld_{23}-\ld_{12})+\ld_{23}(\ld+\ld_3-\ld_{23})+\ld_3(\ld_3+\ld_2+1)+\ld_1^2-\ld\ld_2}\times\\
    \times \dd{\ld_1,\ld_{23},\ld}\dd{\ld_2,\ld_3,\ld_{23}}\dd{\ld_{12},\ld_3,\ld}\dd{\ld_1,\ld_2,\ld_{12}}\prod\limits_{j=0}^{\ld_1+\ld_{23}-\ld-1}\frac{\q{2\ld+2+j}}{\q{2\ld_1-j}\q{2\ld_{23}-j}}\times\\
    \times \frac{1}{\prod\limits_{j=0}^{\ld_{12}+\ld_3-\ld-1}\q{2\ld_{12}-j}}\sum\limits_{s=\ld-\ld_{12}}^{\ld_3}\sum\limits_{z=0}^{\ld_1+\ld_2-\ld_{12}}(-1)^{z}q^{-z(\ld_1+\ld_3+\ld_{12}-\ld_{23})}q^{-s(\ld_3+\ld+\ld_{12}+1)}\times\\
    \times \frac{\q{\ld_2+\ld_3+\ld_{23}+1|\ld_1+\ld_2+\ld_3+2\ld_{23}-\ld_{12}-z-s+1} \q{\ld_2-\ld_1+\ld_{12}|\ld_2+\ld_{23}-z-s}}{\q{z}!\q{\ld_1+\ld_2-\ld_{12}-z}!\q{\ld_1+\ld_{23}-\ld-z}!\q{\ld_{12}-\ld_1-\ld_{23}+z+s}!\q{\ld_3-s}!}\times\\
    \times \frac{\q{\ld_1-\ld_2+\ld_{12}|\ld+\ld_{12}-\ld_2-\ld_3-\ld_{23}+z+s}}{\prod\limits_{j=0}^{\ld_3-\ld_1-\ld_{23}+\ld_{12}+z-1}\q{2\ld_3-j}\prod\limits_{j=0}^{\ld_2-\ld_{23}+s-1}\q{2\ld_2-j}}.
\end{eqnarray*}

\subsubsection{Step 4}
This step is similar to {\it Step 1}. Denote $s=u-z$ and consider summation over $z$. We apply Kronecker delta $\delta_{t,z}$ and summation over all relevant values of $t$, make changes, which we emphasize in {\color{red} red}, and apply identity (\ref{qbid1}).
\begin{eqnarray*}
    &\sum\limits_{z}(-1)^z \frac{q^{-z(\ld_1-\ld_{23}-\ld-1)}}{\q{z}!\q{\ld_1+\ld_2-\ld_{12}-z}!\q{\ld_1+\ld_{23}-\ld-z}!\q{\ld_3-u+z}!\prod\limits_{j=0}^{\ld_3-\ld_1-\ld_{23}+\ld_{12}+z-1}\q{2\ld_3-j}\prod\limits_{j=0}^{\ld_2-\ld_{23}+u-z-1}\q{2\ld_2-j}}=\\
    &=\sum\limits_{z,t,r}(-1)^{z+t+r}\frac{q^{r(t-z)-(r-z)-z(\ld_1-\ld_{23}-\ld-1)-{\color{red} (z-t)(\ld_3-\ld_2-\ld_1+\ld+1)}}}{\q{t-r}!\q{r-z}!\q{z}!\q{\ld_1+\ld_2-\ld_{12}-{\color{red} t}}!\q{\ld_1+\ld_{23}-\ld-{\color{red} t}}!\q{\ld_3-u+z}!}\times\\
    &\times\frac{1}{\prod\limits_{j=0}^{\ld_3-\ld_1-\ld_{23}+\ld_{12}+{\color{red} t}-1}\q{2\ld_3-j}\prod\limits_{j=0}^{\ld_2-\ld_{23}+u-z-1}\q{2\ld_2-j}}.
\end{eqnarray*}
Note that summation over $r,t$ is finite, since $t$ can have values at most $0,1,\ldots,\ld_1+\ld_2-\ld_{12}$ and $r=z,\dots,t$.

\subsubsection{Step 5}
This step is similar to {\it Step 3}. We get rid of the summation over $t$.
\begin{eqnarray*}
    &\sum\limits_{t}(-1)^{t+r}\frac{q^{t(\ld_3-\ld_2-\ld_1+\ld+r+1)}}{\q{t-r}!\q{\ld_1+\ld_2-\ld_{12}-t}!\q{\ld_1+\ld_{23}-\ld-t}!\prod\limits_{j=0}^{\ld_3-\ld_1-\ld_{23}+\ld_{12}+t-1}\q{2\ld_3-j}}=\\
    &=q^{r(\ld_3-\ld_2-\ld_1+\ld+r+1)}q^{-(\ld_2+\ld_1-\ld_{12}-r)(\ld_1+\ld_{23}-\ld-r)}\frac{\prod\limits_{j=0}^{\ld_1+\ld_2-\ld_{12}-r-1}\q{\ld_3+\ld-\ld_{12}-j}}{\q{\ld_1+\ld_2-\ld_{12}-r}!\q{\ld_1+\ld_{23}-\ld-r}!\prod\limits_{j=0}^{\ld_2+\ld_3-\ld_{23}-1}\q{2\ld_3-j}},
\end{eqnarray*}
where we have denoted $t-r=p$ and used Lemma \ref{lemstep5}.

\subsubsection{Step 6}
This step is also similar to {\it Step 3}. We get rid of the summation over $z$.
\begin{eqnarray*}
    \sum\limits_{z}(-1)^z\frac{q^{z(\ld_2-\ld_3+\ld_{23}-r+1)}}{\q{z}!\q{r-z}!\q{\ld_3-u+z}!\prod\limits_{j=0}^{\ld_2-\ld_{23}+u-z-1}\q{2\ld_2-j}}=\\
    =(-1)^rq^{r(\ld_2-\ld_3+\ld_{23}-r+1)}q^{r(\ld_3-u+r)}\frac{\prod\limits_{j=0}^{r-1}\q{\ld_2+\ld_{23}-\ld_3-j}}{\q{r}!\q{\ld_3-u+r}!\prod\limits_{j=0}^{\ld_2-\ld_{23}+u-1}\q{2\ld_2-j}},
\end{eqnarray*}
where we have denoted $r-z=p$ and used Lemma \ref{lemstep6}.

\subsubsection{Intermediate result after Steps 4-6}
From (\ref{qpfid3}) it follows that
\begin{equation*}
    \frac{\q{\ld_2-\ld_1+\ld_{12}|\ld_2+\ld_{23}-u}}{\prod\limits_{j=0}^{\ld_2-\ld_{23}+u-1}\q{2\ld_2-j}}=\q{\ld_2-\ld_1+\ld_{12}|2\ld_2}=\frac{1}{\prod\limits_{j=0}^{\ld_1+\ld_2-\ld_{12}-1}\q{2\ld_2-j}},
\end{equation*}
so
\begin{eqnarray*}
    &\qs{\ld_1}{\ld_2}{\ld_{12}}{\ld_3}{\ld}{\ld_{23}}=(-1)^{\ld_1+\ld_2-\ld_{12}}q^{\ld_3(\ld_3+\ld_{23}+\ld_2+1)+(\ld_1+\ld_{23}-\ld_{12})(\ld+\ld_{12}-\ld_{23}-\ld_2)}\times\\
    &\times \dd{\ld_1,\ld_{23},\ld}\dd{\ld_2,\ld_3,\ld_{23}}\dd{\ld_{12},\ld_3,\ld}\dd{\ld_1,\ld_2,\ld_{12}}\prod\limits_{j=0}^{\ld_1+\ld_{23}-\ld-1}\frac{\q{2\ld+2+j}}{\q{2\ld_1-j}\q{2\ld_{23}-j}}\times\\
    &\times \frac{1}{\prod\limits_{j=0}^{\ld_{12}+\ld_3-\ld-1}\q{2\ld_{12}-j}}\frac{1}{\prod\limits_{j=0}^{\ld_2+\ld_3-\ld_{23}-1}\q{2\ld_3-j}}\frac{1}{\prod\limits_{j=0}^{\ld_1+\ld_2-\ld_{12}-1}\q{2\ld_2-j}}\times\\
    &\times\sum\limits_{u=\ld-\ld_{12}}^{\ld_1+\ld_2+\ld_3-\ld_{12}}\sum\limits_{r=0}^{\ld_1+\ld_2-\ld_{12}}(-1)^{r}q^{r(\ld_1+\ld_2+\ld_3+2\ld_{23}-\ld_{12}+1)}q^{-u(\ld_3+\ld_{12}+\ld+r+1)}\times\\
    &\times \frac{\q{\ld_2+\ld_3+\ld_{23}+1|\ld_1+\ld_2+\ld_3+2\ld_{23}-\ld_{12}-u+1}\q{\ld_1-\ld_2+\ld_{12}|\ld+\ld_{12}-\ld_2-\ld_3-\ld_{23}+u}}{\q{\ld_{12}-\ld_1-\ld_{23}+u}!\q{\ld_3-u+r}!}\\
    &\times \frac{\prod\limits_{j=0}^{\ld_1+\ld_2-\ld_{12}-r-1}\q{\ld_3+\ld-\ld_{12}-j}\prod\limits_{j=0}^{r-1}\q{\ld_2+\ld_{23}-\ld_3-j}}{\q{r}!\q{\ld_1+\ld_2-\ld_{12}-r}!\q{\ld_1+\ld_{23}-\ld-r}!}.
\end{eqnarray*}

\subsubsection{Step 7} In this step we will get rid of the summation over $u$. Consider
\begin{eqnarray*}
    &\sum\limits_{u}q^{-u(\ld_3+\ld_{12}+\ld+r+1)}\frac{\q{\ld_2+\ld_3+\ld_{23}+1|\ld_1+\ld_2+\ld_3+2\ld_{23}-\ld_{12}-u+1}\q{\ld_1-\ld_2+\ld_{12}|\ld+\ld_{12}-\ld_2-\ld_3-\ld_{23}+u}}{\q{\ld_{12}-\ld_1-\ld_{23}+u}!\q{\ld_3-u+r}!}=\\
    &=q^{-(\ld_3+r)(\ld_3+\ld_{12}+\ld+r+1)}q^{(\ld_{12}-\ld_1-\ld_{23}+\ld_3+r)(\ld+\ld_{12}-\ld_2-\ld_{23}+r)}\times\\
    &\times\frac{\q{\ld_1-\ld_2+\ld_{12}|\ld+\ld_{12}-\ld_2-\ld_{23}+r}\prod\limits_{j=0}^{\ld_{12}-\ld_{23}+\ld_3-\ld_1+r-1}\q{\ld_3+\ld+\ld_{12}+r+1-j}}{\q{\ld_{12}-\ld_1-\ld_{23}+\ld_3+r}!}
\end{eqnarray*}

\subsubsection{The result} 
Taking into account that
\begin{eqnarray*}
    &\prod\limits_{j=0}^{\ld_1+\ld_{23}-\ld-1}\q{2\ld+2+j} \prod\limits_{j=0}^{\ld_{12}-\ld_{23}+\ld_3-\ld_1+r-1}\q{\ld_3+\ld+\ld_{12}+r+1-j}=\prod\limits_{j=0}^{\ld_3+\ld_{12}-\ld+r-1}\q{2\ld+2+j},
\end{eqnarray*}
\begin{eqnarray*}
    \frac{\q{\ld_1-\ld_2+\ld_{12}|\ld+\ld_{12}-\ld_2-\ld_{23}+r}}{\prod\limits_{j=0}^{\ld_1+\ld_{23}-\ld-1}\q{2\ld_1-j}}=\frac{\prod\limits_{j=0}^{\ld_1+\ld_2-\ld_{12}-r-1}\q{\ld_1-\ld_{23}+\ld-j}}{\prod\limits_{j=0}^{\ld_1+\ld_2-\ld_{12}-1}\q{2\ld_1-j}}
\end{eqnarray*}
and substituting $\ld_1+\ld_2-\ld_{12}-r=z$, we obtain
\begin{eqnarray*}
    &\qs{\ld_1}{\ld_2}{\ld_{12}}{\ld_3}{\ld}{\ld_{23}}= \frac{\dd{\ld_1,\ld_{23},\ld}}{\prod\limits_{j=0}^{\ld_1+\ld_{23}-\ld-1}\q{2\ld_{23}-j}}\frac{\dd{\ld_2,\ld_3,\ld_{23}}}{\prod\limits_{j=0}^{\ld_2+\ld_3-\ld_{23}-1}\q{2\ld_3-j}}\frac{\dd{\ld_{12},\ld_3,\ld}}{\prod\limits_{j=0}^{\ld_{12}+\ld_3-\ld-1}\q{2\ld_{12}-j}}\frac{\dd{\ld_1,\ld_2,\ld_{12}}}{\prod\limits_{j=0}^{\ld_1+\ld_2-\ld_{12}-1}\q{2\ld_1-j}}\times\nonumber\\
    &\times \frac{1}{\prod\limits_{j=0}^{\ld_1+\ld_2-\ld_{12}-1}\q{2\ld_2-j}}\sum\limits_{z=0}^{\ld_1+\ld_2-\ld_{12}}(-1)^z \prod\limits_{j=0}^{\ld_1+\ld_2+\ld_3-\ld-z-1}\q{2\ld+2+j}\times\\
    &\times\frac{\prod\limits_{j=0}^{z-1}\q{\ld+\ld_1-\ld_{23}-j}\prod\limits_{j=0}^{z-1}\q{\ld-\ld_{12}+\ld_3-j}\prod\limits_{j=0}^{\ld_1+\ld_2-\ld_{12}-z-1}\q{\ld_2+\ld_{23}-\ld_3-j}}{\q{z}!\q{\ld_1+\ld_2-\ld_{12}-z}!\q{\ld_2+\ld_3-\ld_{23}-z}!\q{\ld_{12}+\ld_{23}-\ld_2-\ld+z}!}\nonumber.
\end{eqnarray*}

\section{Properties of quantum $6j$-symbols}\label{sectionpropq6j}
From now on for any two finite disjoint sets $I$, $I^\prime$ consisting of natural numbers and their union $I\cup I^\prime$, we assume that $\ld_I$, $\ld_{I^\prime}, \ld_{I\cup I^\prime}$ are generic and $\ld_I+\ld_{I^\prime}-\ld_{I\cup I^\prime}\in\mbb{Z}_{\geq 0}$. We also denote
$$
C_I:=\ld_I(\ld_I+1).
$$
\subsection{Formula for $q6j$-symbols}
\begin{theorem}
$q6j$-symbols in (\ref{qsdef}) are given by formula
\begin{eqnarray}\label{q6jfinal}
    &\qs{\ld_1}{\ld_2}{\ld_{12}}{\ld_3}{\ld_{123}}{\ld_{23}}= \frac{\dd{\ld_1,\ld_{23},\ld_{123}}}{\prod\limits_{j=0}^{\ld_1+\ld_{23}-\ld_{123}-1}\q{2\ld_{23}-j}}\frac{\dd{\ld_2,\ld_3,\ld_{23}}}{\prod\limits_{j=0}^{\ld_2+\ld_3-\ld_{23}-1}\q{2\ld_3-j}}\frac{\dd{\ld_{12},\ld_3,\ld_{123}}}{\prod\limits_{j=0}^{\ld_{12}+\ld_3-\ld_{123}-1}\q{2\ld_{12}-j}}\frac{\dd{\ld_1,\ld_2,\ld_{12}}}{\prod\limits_{j=0}^{\ld_1+\ld_2-\ld_{12}-1}\q{2\ld_1-j}}\times\nonumber\\
    &\times \frac{1}{\prod\limits_{j=0}^{\ld_1+\ld_2-\ld_{12}-1}\q{2\ld_2-j}}\sum\limits_{z=0}^{\ld_1+\ld_2-\ld_{12}}(-1)^z \prod\limits_{j=0}^{\ld_1+\ld_2+\ld_3-\ld_{123}-z-1}\q{2\ld_{123}+2+j}\times\\
    &\times\frac{\prod\limits_{j=0}^{z-1}\q{\ld_{123}+\ld_1-\ld_{23}-j}\prod\limits_{j=0}^{z-1}\q{\ld_{123}-\ld_{12}+\ld_3-j}\prod\limits_{j=0}^{\ld_1+\ld_2-\ld_{12}-z-1}\q{\ld_2+\ld_{23}-\ld_3-j}}{\q{z}!\q{\ld_1+\ld_2-\ld_{12}-z}!\q{\ld_2+\ld_3-\ld_{23}-z}!\q{\ld_{12}+\ld_{23}-\ld_2-\ld_{123}+z}!}\nonumber.
\end{eqnarray}
\end{theorem}
\begin{proof}
    Expression (\ref{q6jfinal}) has been obtained as a result of Section \ref{deriveformulaq6j}.
\end{proof}

\subsection{The shadow world}\label{shadworld}
\begin{figure}[h]
\begin{center}
	{\includegraphics[width=260pt]{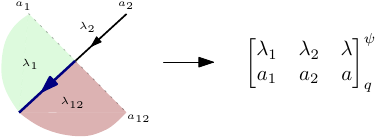}}
	\caption{Graphical representation of $q3j$-symbols in the shadow world.}
	\label{crossing0}
\end{center}
\end{figure}

In the following subsections we will examine properties of $q6j$-symbols, proofs for which rely on transitioning to the shadow world. Lemma \ref{shadowtransit} is the key ingredient for this transition.
\begin{lemma}\label{shadowtransit}
\begin{eqnarray}\label{stid1}
    \sum\limits_{a_1,a_2,a_3}\qt{\ld_2}{\ld_3}{\ld_{23}}{a_2}{a_3}{a_{23}}^\pi\qt{\ld_{12}}{\ld_3}{\ld_{123}}{a_{12}}{a_3}{a_{123}}^\psi\qt{\ld_1}{\ld_2}{\ld_{12}}{a_1}{a_2}{a_{12}}^\psi=\\
    =\qs{\ld_1}{\ld_2}{\ld_{12}}{\ld_3}{\ld_{123}}{\ld_{23}}\qt{\ld_1}{\ld_{23}}{\ld_{123}}{a_1}{a_{23}}{a_{123}}^\psi\nonumber,
\end{eqnarray}
\begin{figure}[h]
\begin{center}
	{\includegraphics[width=320pt]{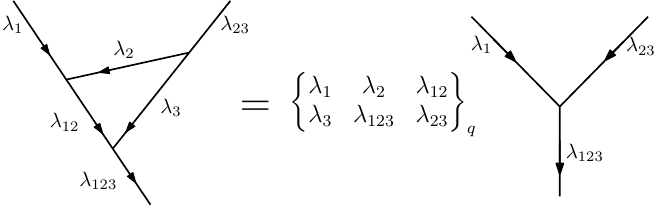}}
    \caption{String diagram representation for the relation (\ref{stid1}).}
\end{center}
\end{figure}\par
\begin{eqnarray}\label{stid2}
\sum\limits_{a_1,a_2,a_3}\qt{\ld_1}{\ld_{23}}{\ld_{123}}{a_1}{a_{23}}{a_{123}}^\psi\qt{\ld_2}{\ld_3}{\ld_{23}}{a_2}{a_3}{a_{23}}^\psi=\sum\limits_{\ld_{12}}\qs{\ld_3}{\ld_2}{\ld_{23}}{\ld_1}{\ld_{123}}{\ld_{12}}\times\\
\times \sum\limits^\prime_{a_1,a_2,a_3}\qt{\ld_{12}}{\ld_3}{\ld_{123}}{a_{12}}{a_3}{a_{123}}^\psi\qt{\ld_1}{\ld_2}{\ld_{12}}{a_1}{a_2}{a_{12}}^\psi\nonumber,
\end{eqnarray}
\begin{figure}[h]
\begin{center}
	{\includegraphics[width=330pt]{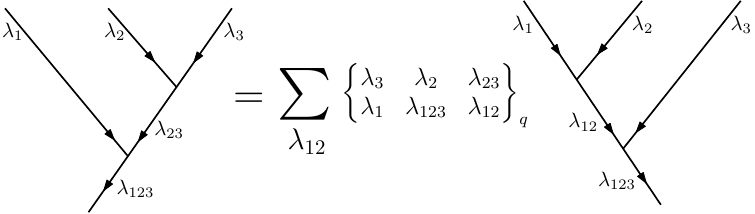}}
    \caption{String diagram representation for the relation (\ref{stid2}).}
\end{center}
\end{figure}\par
\begin{eqnarray}\label{stid3}
\sum\limits_{a_1,a_2,a_3}\sum\limits_{a_2^\prime,a_3^\prime}\left(\mathcal{R}^{\ld_2,\ld_3}\right)^{a_2^\prime,a_3^\prime}_{a_2,a_3}\qt{\ld_{13}}{\ld_2}{\ld_{123}}{a_{13}}{a^\prime_2}{a_{123}}^\psi\qt{\ld_1}{\ld_3}{\ld_{13}}{a_1}{a^\prime_3}{a_{13}}^\psi=\\
=\sum\limits_{\ld_{12}}(-1)^{\ld_{13}+\ld_{12}-\ld_{123}-\ld_1}q^{C_{123}+C_{1}-C_{13}-C_{12}}\times\nonumber\\
\times\qs{\ld_3}{\ld_1}{\ld_{13}}{\ld_2}{\ld_{123}}{\ld_{12}}\sum\limits^\prime_{a_1,a_2,a_3}\qt{\ld_{12}}{\ld_3}{\ld_{123}}{a_{12}}{a_3}{a_{123}}^\psi\qt{\ld_1}{\ld_2}{\ld_{12}}{a_1}{a_2}{a_{12}}^\psi\nonumber.
\end{eqnarray}
\begin{figure}[h]
\begin{center}
	{\includegraphics[width=430pt]{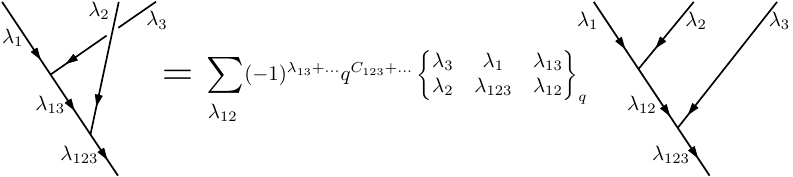}}
    \caption{String diagram representation for the relation (\ref{stid3}).}
    \label{crossing1}
\end{center}
\end{figure}\par
\end{lemma}
\begin{proof}
See Appendix \ref{prooftr}.
\end{proof}
Let us represent $q3j$-symbols as in Figure \ref{crossing0}, where it is enough to consider embeddings $\psi$.
\par
\begin{figure}[h]
\begin{center}
	{\includegraphics[width=330pt]{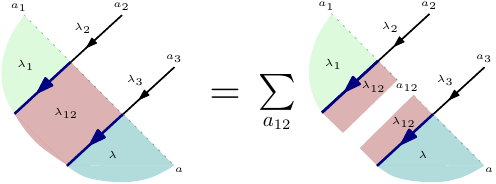}}
	\caption{Graphical representation of the composition of $q3j$-symbols in the shadow world.}
	\label{crossing00}
\end{center}
\end{figure}

The gray dotted line divides the plane in two parts. In the upper part we have string diagrams defined in Section \ref{strings}. The lower part is called \textit{the shadow world}. It consists of strings and gleams between them. Strings do not change their colors and orientations after crossing the dotted line, but their ends are no longer marked. Gleams are colored by highest weights $\ld_j$ of Verma modules $M_{\ld_j}$, in accordance to the colors of the gray dotted lines before the transition. For example, in the second row of Figure \ref{assign} the string colored by $\ld_1$ with its end marked by $a_1$ now became a gray line with its end marked by $a_1$ and with a green gleam underneath, which is colored by $\ld_1$. Applying this assignment to Figure \ref{assign3}, the composition of $q3j$-symbols in the shadow world is represented as in Figure \ref{crossing00}.
\par

\begin{figure}[h]
\begin{center}
	{\includegraphics[width=320pt]{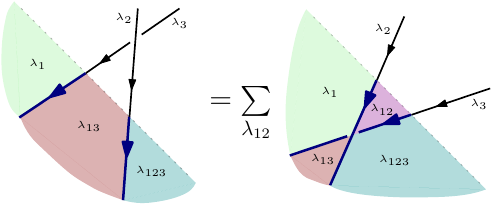}}
	\caption{Graphical representation of the identity (\ref{stid3}) in the shadow world.}
	\label{crossing2}
\end{center}
\end{figure}

\begin{figure}[h]
\begin{center}
	{\includegraphics[width=440pt]{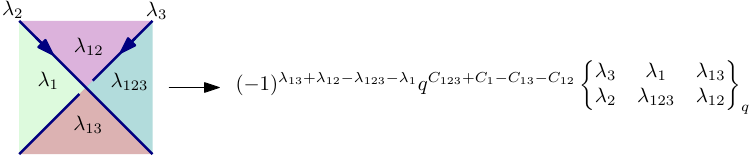}}
	\caption{Weight assignment for the crossing in the shadow world.}
	\label{crossing3}
\end{center}
\end{figure}

Lemma \ref{shadowtransit} is used to determine weights assigned to elementary graphs in the shadow world: (\ref{stid1}) defines trivalent vertex with one end splitting in two, (\ref{stid2}) defines trivalent vertex with two ends joining in one, (\ref{stid3}) defines crossings. As an example we show how this procedure works for the identity (\ref{stid3}), as shown in Figures \ref{crossing1}, \ref{crossing2}, \ref{crossing3}.\par

\par
We use the rules described above to identify weights given in Figure \ref{crossing1} with their graphical representations in the shadow world given in Figure \ref{crossing2}. The result is shown in Figure \ref{crossing3}.
\par
\begin{figure}[h]
\begin{center}
	{\includegraphics[width=460pt]{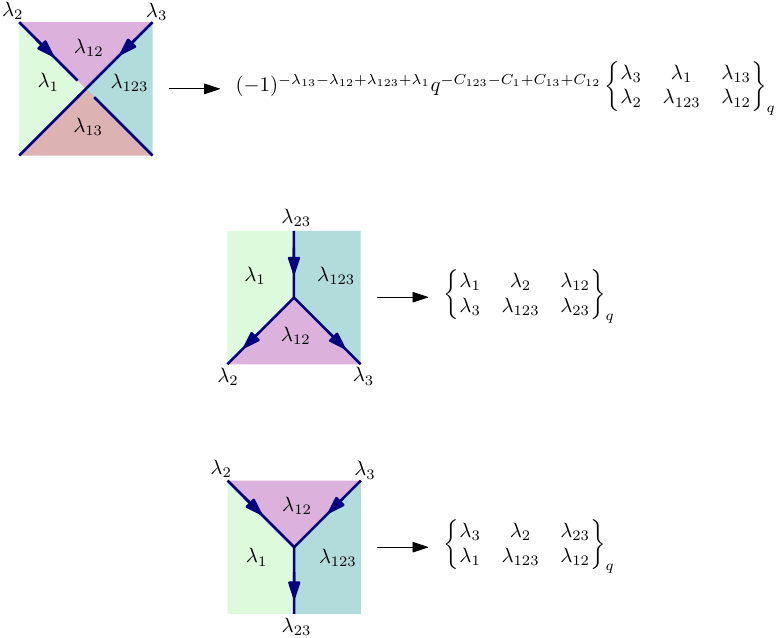}}
    \caption{Elementary graphs in the shadow world and their weights: the first line assigns weights given by identity (\ref{stid3}) for the inverse $\mathcal{R}$-matrix to the opposite crossing, the second line assigns weights given by identities (\ref{stid1}) and (\ref{stid2}) to the corresponding trivalent vertices.}
    \label{trivv}
\end{center}
\end{figure}

One can use the identity (\ref{stid3}) with the inverse $\mathcal{R}$-matrix to obtain weights assigned to the opposite crossings in the shadow world. Similarly, applying the same procedure to the identities (\ref{stid1}) and (\ref{stid2}) yields the weights for the corresponding trivalent vertices. The resulting weight assignments for all elementary graphs in the shadow world are summarized in Figures \ref{crossing3} and \ref{trivv}.
\par
This procedure provides a systematic way to translate any identity from Section \ref{qtprop} into the shadow world framework. In the following subsections, we adopt this approach as our main strategy for deriving explicit formulae and establishing properties of $q6j$-symbols.\par

\subsection{Orthogonality relation}\label{ortsect}
\begin{theorem}
The following identity holds
\begin{equation}\label{qsorth}
    \sum_{\alpha}\qs{\ld_1}{\ld_2}{\alpha}{\ld_3}{\ld_{123}}{\mu}\qs{\ld_3}{\ld_2}{\nu}{\ld_1}{\ld_{123}}{\alpha}=\delta_{\mu,\nu}.
\end{equation}
\end{theorem}
\begin{figure}[h]
\begin{center}
	{\includegraphics[width=280pt]{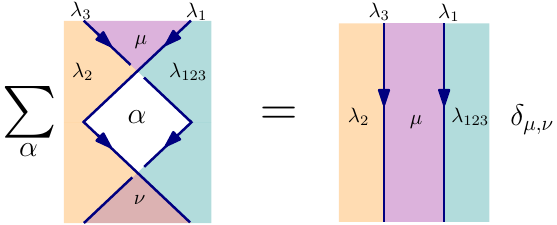}}
    \caption{Graphical representation of the orthogonality relation (\ref{qsorth}) in the shadow world, using the rules from Figures \ref{crossing3} and \ref{trivv}.}
\end{center}
\end{figure}\par
\begin{proof}
    The proof parallels \cite[Theorem 6.1, Eq. (6.16)]{KR}, but uses Lemma \ref{shadowtransit} and the shadow world transition of Theorem \ref{orthrelth}, identity (\ref{projth2}).
\end{proof}
Note that the summation in (\ref{qsorth}) is over $\alpha=\ld_1+\ld_2-j$, where $j=0,1,\ldots,\ld_1+\ld_2+\ld_3-\ld_{123}$.

\subsection{Racah identity}\label{racahsect}
\begin{theorem}
The following identity holds
\begin{eqnarray}\label{qsracah}
    &\sum\limits_\alpha (-1)^\alpha q^{-\alpha(\alpha+1)}\qs{\ld_1}{\ld_2}{\alpha}{\ld_3}{\ld_{123}}{\ld_{23}}\qs{\ld_3}{\ld_1}{\ld_{13}}{\ld_2}{\ld_{123}}{\alpha}=\\
    &=(-1)^{\ld_1+\ld_2+\ld_3+\ld_{123}-\ld_{13}-\ld_{23}} q^{C_{13}+C_{23}-C_1-C_2-C_3-C_{123}} \qs{\ld_1}{\ld_3}{\ld_{13}}{\ld_2}{\ld_{123}}{\ld_{23}}.\nonumber
\end{eqnarray}
\end{theorem}
\begin{figure}[h]
\begin{center}
	{\includegraphics[width=400pt]{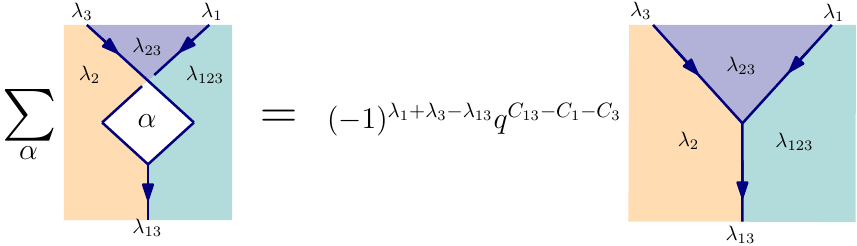}}
    \caption{The shadow world representation of the Racah identity (\ref{qsracah}).}
\end{center}
\end{figure}\par
\begin{proof}
    We refer to the proof of \cite{KR}, Theorem 6.1, formula (6.17), but one needs to use Lemma \ref{shadowtransit} and Theorem \ref{rmatth}, identity (\ref{rmatid1}) from the present paper.
\end{proof}
Note that the summation in (\ref{qsracah}) is also over $\alpha=\ld_1+\ld_2-j$, where $j=0,1,\ldots,\ld_1+\ld_2+\ld_3-\ld_{123}$.

\subsection{Biedenharn-Elliot identity}\label{BEsect}
This identity is also known as pentagon identity.
\begin{theorem}
The following identity holds
\begin{eqnarray}\label{qsBE}
    \sum\limits_\alpha \qs{\ld_1}{\ld_2}{\ld_{12}}{\ld_3}{\ld_{123}}{\alpha}\qs{\ld_1}{\alpha}{\ld_{123}}{\ld_4}{\ld_{1234}}{\ld_{234}}\qs{\ld_2}{\ld_3}{\alpha}{\ld_4}{\ld_{234}}{\ld_{34}}=\\
    =\qs{\ld_{12}}{\ld_3}{\ld_{123}}{\ld_4}{\ld_{1234}}{\ld_{34}}\qs{\ld_1}{\ld_2}{\ld_{12}}{\ld_{34}}{\ld_{1234}}{\ld_{234}}.\nonumber
\end{eqnarray}
\end{theorem}
\begin{figure}[h]
\begin{center}
	{\includegraphics[width=390pt]{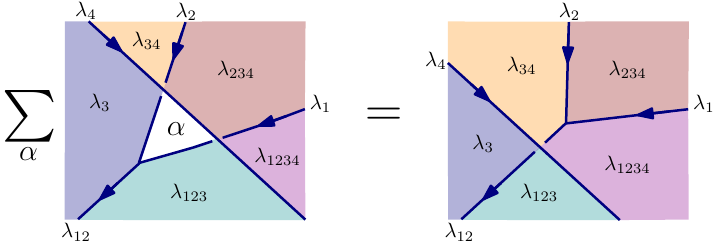}}
    \caption{The shadow world representation of the Biedenharn-Elliot identity (\ref{qsBE}).}
\end{center}
\end{figure}\par
\begin{proof}
    We refer to \cite{KR}, Theorem 6.1, formula (6.18), but one needs to use Lemma \ref{shadowtransit} and Theorem \ref{rmatth}, identity (\ref{rmatid2}).
\end{proof}
Note that the summation in (\ref{qsBE}) is over $\alpha=\ld_2+\ld_3-j$, where $j=0,1,\ldots,\mrm{min}\{\ld_2+\ld_3+\ld_4-\ld_{234},\ld_1+\ld_2+\ld_3-\ld_{123}\}$. Also note that $\ld_{12}+\ld_{34}-\ld_{1234}\geq 0$.

\subsection{Yang-Baxter identity}\label{YBsect}
\begin{theorem}
The following identity holds
\begin{eqnarray}\label{qsYB}
    &\sum\limits_{\alpha}(-1)^{-\alpha+\ld_2-\ld_{23}-\ld_{24}}q^{-\alpha(\alpha+1)+C_2-C_{23}-C_{24}}\times\\
    &\times \qs{\ld_1}{\ld_2}{\alpha}{\ld_3}{\ld_{123}}{\ld_{23}}\qs{\ld_4}{\alpha}{\ld_{124}}{\ld_3}{\ld_{1234}}{\ld_{123}}\qs{\ld_4}{\ld_2}{\ld_{24}}{\ld_1}{\ld_{124}}{\alpha}=\nonumber\\
    &=\sum\limits_{\beta}(-1)^{-\beta+\ld_{1234}-\ld_{123}-\ld_{124}}q^{-\beta(\beta+1)+C_{1234}-C_{123}-C_{124}}\times\nonumber\\
    &\times \qs{\ld_4}{\ld_{23}}{\beta}{\ld_1}{\ld_{1234}}{\ld_{123}}\qs{\ld_4}{\ld_2}{\ld_{24}}{\ld_3}{\beta}{\ld_{23}}\qs{\ld_1}{\ld_{24}}{\ld_{124}}{\ld_3}{\ld_{1234}}{\beta}.\nonumber
\end{eqnarray}
\end{theorem}
\begin{figure}[h]
\begin{center}
	{\includegraphics[width=380pt]{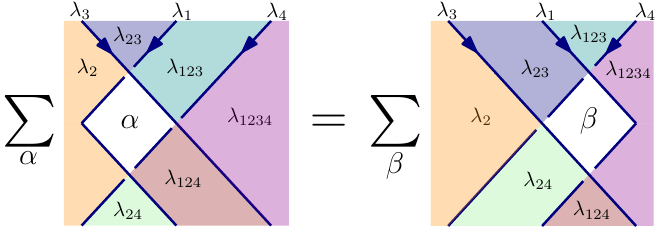}}
    \caption{The shadow world representation of the Yang-Baxter identity (\ref{qsYB}).}
\end{center}
\end{figure}\par
\begin{proof}
    We refer to \cite{KR}, Theorem 6.1, formula (6.19), but one needs to use Lemma \ref{shadowtransit} and (\ref{YBR}).
\end{proof}
Note that in (\ref{qsYB}) the summation in the left hand side is over $\alpha=\ld_1+\ld_2-j$, $j=0,\ldots,\mrm{min}\{\ld_1+\ld_2+\ld_3-\ld_{123},\ld_1+\ld_2+\ld_4-\ld_{124}\}$, and in the right hand side is over $\beta=\ld_2+\ld_3+\ld_4-k$, $k=0,1,\ldots,\ld_1+\ld_2+\ld_3+\ld_4-\ld_{1234}$.

Properties presented above can also be proven using ideas explored in \cite{T}, Chapter VI, Section 5. They rely on multiplicity modules and their tensor contractions. See Theorems 5.4.2, 5.4.3, 5.4.1, 5.7 respectively to Sections \ref{ortsect}, \ref{racahsect}, \ref{BEsect}, \ref{YBsect}.

Similarly to $q3j$-symbols, with proper restrictions being imposed, the obtained formulae for $q6j$-symbols can be applied to tensor products with components being finite-dimensional irreducible representations, given that highest weights of the components are chosen generically and the decomposition consists of irreducible representations only. Note that $q6j$-symbols for Verma modules admit fewer symmetries, i.e. there is no Regge symmetry since such tensor product decomposition wouldn't make sense. We also didn't explore analogs of identities (6.10) and (6.11) in \cite{KR}, which compute weights of cups and caps, since there are no corresponding identities to transpose to the shadow world. This issue is discussed in detail in Section \ref{chevv}.
\section{On $F_K(x,q)$ and the shadow world}\label{swsec}
In this section, we will use the shadow world graphical calculus to outline a construction of a functional $\mathsf{SW}_\infty$, which is a possible model for Gukov-Manolescu series $F_K(x,q)$. This section contains preliminary results, which will be expanded upon in \cite{S}.

\subsection{Reverse-engineering evaluation and coevaluation maps}\label{chevv}
\subsubsection{Outline of the problem}
In Definition \ref{defshap} we used an involution (\ref{rho}) to remain within the $q$-analog of the BGG category $\mathcal{O}$, which allowed us to do the normalization procedure for $q3j$-symbols using only Verma modules of the highest weight. However, the Shapovalov form defined this way cannot be used as a $\Uq$-invariant evaluation map, and hence it cannot be represented as a string diagram in a way consistent with the topology of graphs. Similarly, $\Uq$-invariance of the coevaluation map would imply going beyond Verma modules of the highest weight.\par

\begin{figure}[h]
\begin{center}
	{\includegraphics[width=340pt]{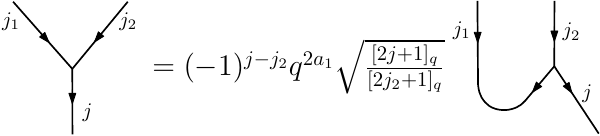}}
	\caption{An identity relating an embedding $\psi: V_j\to V_{j_1}\ot V_{j_2}$ with a composition of projection $\pi: V_{j_1}\ot V_{j} \to V_{j_2}$ and coevaluation map $i:\mathbb{C}\to V_{j_1}\ot V_{j_1}$.}
	\label{bending}
\end{center}
\end{figure}

To illustrate this problem, consider a tensor product of two finite-dimensional modules $V_{j_1}\ot V_{j_2}$, where $j_1,j_2\in\frac{1}{2}\mbb{Z}_{\geq 0}$. There is an identity, which relates an embedding $\psi: V_j\to V_{j_1}\ot V_{j_2}$ with a projection $\pi: V_{j_1}\ot V_{j} \to V_{j_2}$ and the coevaluation map $i:\mathbb{C}\to V_{j_1}\ot V_{j_1}$, which corresponds to the topological move of string diagrams as in Figure \ref{bending}.
\par
For the case of tensoring two Verma modules of the highest weight $M_{\ld_1}\ot M_{\ld_2}$ and given that $\ld=\ld_1+\ld_2-J$ for some $J\in\mbb{Z}_{\geq 0}$ the right hand side of the identity above makes no sense. Formula (\ref{prj}) for projections cannot be used for such values of arguments. From the representation theory perspective, there is no such projection $\pi:M_{\ld_1}\ot M_{\ld} \to M_{\ld_2}$ since there is no component $M_{\ld_2}$ in the tensor product decomposition of $M_{\ld_1}\ot M_{\ld}$.\par

For a similar reason, we cannot express the coevaluation map and the evaluation map using embeddings and projections, respectively, by putting $j_1=-j_2$, $j=0$, since this would involve both highest weight and lowest weight Verma modules.\par

Inability to use the identity in Figure \ref{bending} is an obstacle because it is a key ingredient for definition of cups and caps in the shadow world. As a result of Section \ref{shadworld}, we managed to obtain graphical assignment only for crossings and trivalent vertices, as in Figures \ref{crossing3} and \ref{trivv}.

\subsubsection{The Chevalley extended $\widetilde{\Uq}$ as a solution}\label{chev}
\begin{figure}[h]
\begin{center}
	{\includegraphics[width=500pt]{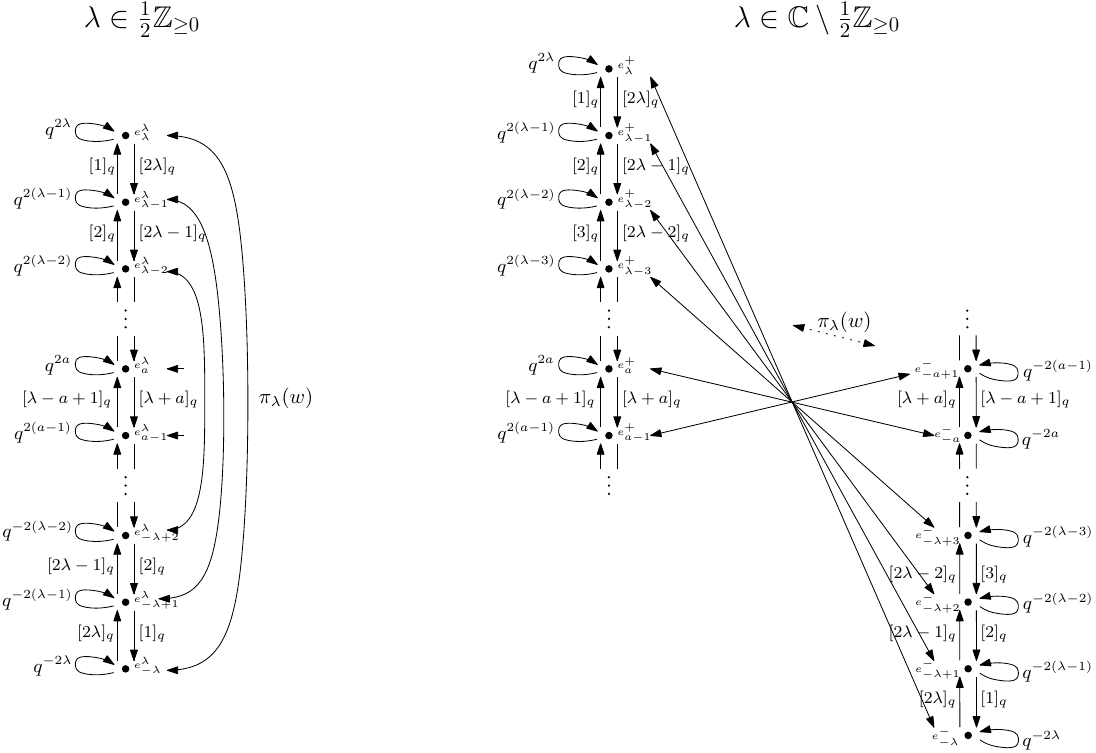}}
	\caption{Weight diagram of a polarized module $P_\lambda$. For $\ld\in\mbb{C}\setminus \frac{1}{2}\Zp$ it is an incarnation of a large color limit of a finite-dimensional $\Uq$-module.}
	\label{polmod}
\end{center}
\end{figure}

Let us briefly sketch the solution for the problem outlined above. To proceed further, it is necessary to use $\Uq$-invariant evaluation $(-,-)_\pm$ and coevaluation $i_\pm$ maps, which are defined for pairs of Verma module $M_\ld^+$ of the highest weight $\ld$ with Verma module $M_{-\ld}^-$ of the lowest weight $-\ld$. These two modules are related by the Chevalley involution, which can be made into an inner automorphism by extending $\Uq$ to a $\widetilde{\Uq}$, where we add a $q$-analog of the Weyl element $w$. For a proper account of this construction, see, for example, \cite{R2}. Hence, we obtain a $\widetilde{\Uq}$-module $P_\ld$, such that 
\begin{equation}\label{polprop}
P_\ld \big|_{\Uq}=
\begin{cases}
M^+_\ld\oplus M^-_{-\ld},& \ld\in\mbb{C}\setminus \frac{1}{2}\Zp,\\
V_\ld,& \ld\in\frac{1}{2}\Zp.
\end{cases}
\end{equation}
We will refer to this module as a \textit{polarized module}. For $\ld\in \mbb{C}\setminus \frac{1}{2}\Zp$, the polarized module $P_\ld$ can be understood as a Weyl-symmetric analytic continuation of finite-dimensional $\Uq$-modules. In particular, it extends many of their structural features while remaining well-defined for non-integral highest weights. The corresponding weight diagram is presented in Figure \ref{polmod}.
\par

\begin{figure}[h]
\begin{center}
	{\includegraphics[width=460pt]{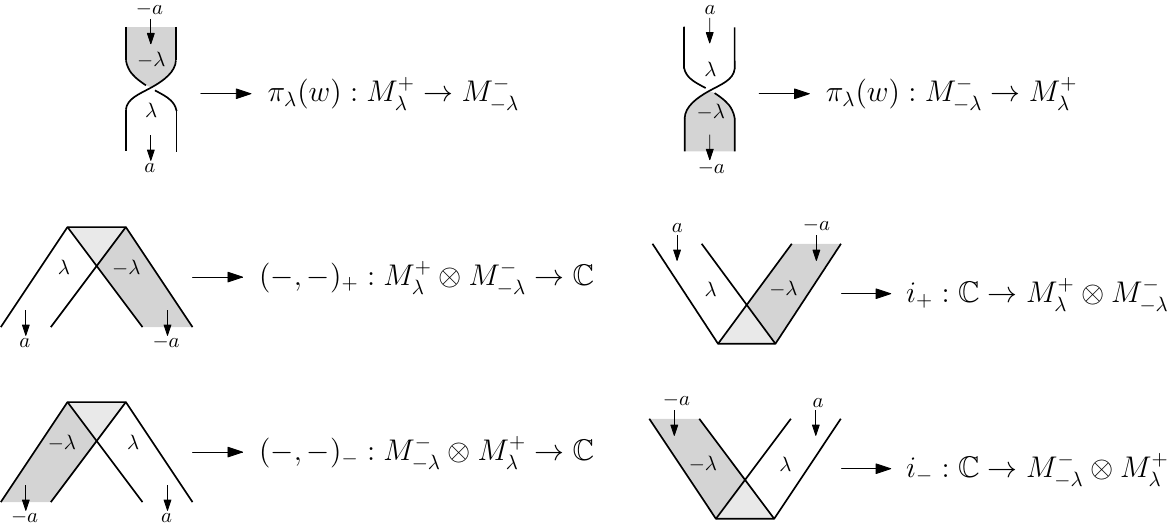}}
	\caption{The first row shows ribbon graphs representing the Weyl element $w$ in the representation $P_\ld$ of $\widetilde{\Uq}$, realized through the representation map $\pi_\ld$. In this setting, the highest weight Verma module $M^+_\ld$ and the lowest weight Verma module $M^-_{-\ld}$ appear as $\Uq$-submodules of $P_\ld$, connected to each other by the action of $\pi_\ld(w)$. The rows below illustrate the ribbon graphs corresponding to the evaluation and coevaluation maps.}
	\label{weyl1}
\end{center}
\end{figure}

This motivates upgrading the string diagrams introduced in Section \ref{strings} to ribbon graphs. The ribbons come in two types: white, representing highest weight Verma modules, and gray, representing lowest weight Verma modules. These two are connected by a twist, which graphically encodes the action of the Weyl element. With this framework, the evaluation and coevaluation maps can be represented as shown in Figure \ref{weyl1}.
\begin{remark}\label{remarkpm}
Since $\Uq$-submodules $M^+_\ld$ and $M^-_{-\ld}$ of a $\widetilde{\Uq}$-module $P_\ld$ are not related by operators $E$ and $F$, there should be a definition of two different evaluation and coevaluation maps, i.e. $i_+:\mathbb{C}\to M^+_\ld \ot M^-_{-\ld}$ and $i_-:\mathbb{C}\to M^-_{-\ld} \ot M^+_{\ld}$.
\end{remark}
\par

\begin{figure}[h]
\begin{center}
	{\includegraphics[width=290pt]{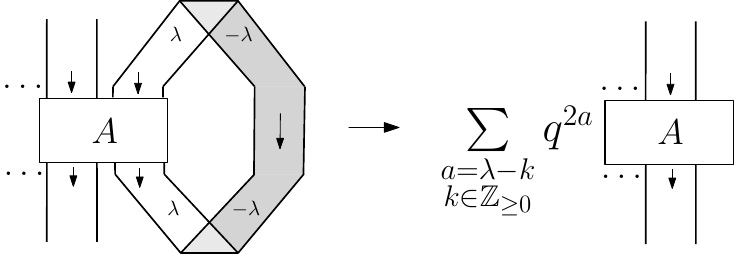}}
	\caption{Ribbon diagram for taking a trace of operator $A$ of a general form over a single ribbon strand colored by a highest weight $\ld$, where "$\ldots$" stands for other strands.}
	\label{weyl2}
\end{center}
\end{figure}

Once taking a trace of a braid colored by Verma modules of the highest weight, one would need to repeatedly apply a composition of evaluation and coevaluation maps $(-,-)_+$ and $i_+$, as shown in Figure \ref{weyl2}, where as an example we have taken a trace over a single ribbon strand colored by $M^+_\ld$. This amounts to adding a $q^{2(\ld-k)}$ factor to the matrix element of the operator $A$ and summing over one of its indices $k\in\mathbb{Z}_{\geq 0}$.
\par
This shows that if we restrict ourselves to braids colored by Verma modules of the highest weight, we can still take traces in a functorial way without introducing ribbon strands colored by Verma modules of the lowest weight. A string diagram of the junction of a cup and a cap corresponds to a factor $q^{2a}$ and a summation over the weight space of $M^+_{\ld}$. However, this assignment is ad hoc and its usage is limited to taking traces of braids, while certain topological moves remain unavailable, i.e. rotating the $\mathcal{R}$-matrix and bending strings of $q3j$-symbols.\par

Since Verma modules of the highest weight were the main focus of previous sections, in what follows we will be considering them first. Nevertheless, for the purpose of reproducing the GM series $F_K(x,q)$, it is more natural to work with braids colored by polarized modules $P_\ld$, $\ld\in \mbb{C}\setminus \frac{1}{2}\Zp$. A possible way to proceed with this is to also consider the same braid colored by the lowest weight Verma modules $M^-_{-\ld}$, where we have "inverted" the initial colors $M^+_\ld$, and take a trace in a similar fashion using a junction of cups and caps given by $i_-$ and $(-,-)_-$ respectively, instead of $i_+$ and $(-,-)_+$ used previously. For each strand, such a junction contributes a factor $q^{-2(\ld-k)}$ to the matrix element and a summation over $k\in\Zp$, which amounts to a summation over the weight space of $M^-_{-\ld}$. The sum of traces of braids colored by modules $M^+_\ld$ and their “inverse colors” $M^-_{-\ld}$ amounts to taking a trace of a braid colored by the corresponding polarized modules $P_\ld$.\par

\subsubsection{Reverse-engineering cups and caps in the shadow world}
\begin{figure}[h]
\begin{center}
	{\includegraphics[width=420pt]{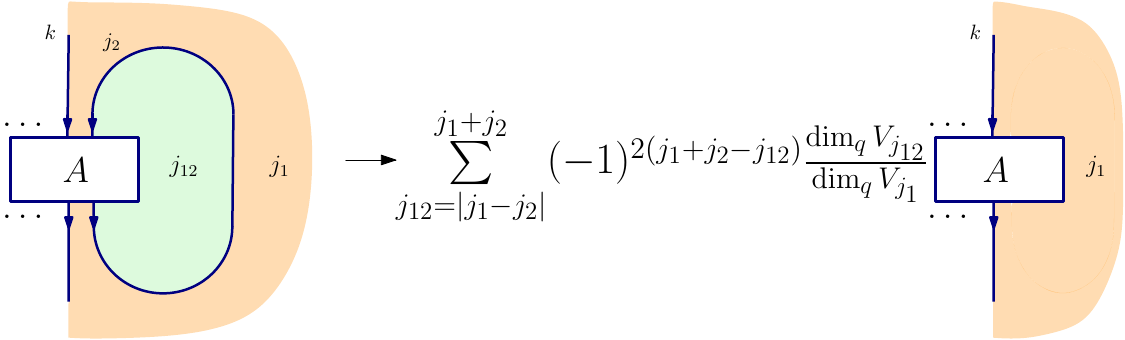}}
	\caption{Strands are colored by finite-dimensional representations $V_k$ and $V_{j_2}$, while gleams are colored by $V_{j_1}$ and $V_{j_{12}}$. In the right hand side we have weights assigned to a composition of a cup and a cap in the shadow world.}
	\label{weylSW}
\end{center}
\end{figure}

Consider the diagram depicted in Figure \ref{weyl2}, where the trace is taken over a strand colored by a finite-dimensional representation $V_j$, $j\in\frac{1}{2}\Zp$. In this case, the contributing factor is also given by $q^{2a}$, but the summation will be over the weight space of $V_j$ instead of $M^+_\ld$. In the shadow world, we have a similar diagram, as depicted in Figure \ref{weylSW}.
\\
Here
$$
\mathrm{dim}_q \,V_j:=\mrm{Tr}_{V_j}(K)=\sum\limits_{a=-j}^j q^{2a}=\q{2j+1}.
$$
In this figure we have taken into account the cup and the cap for the strand colored by $V_{j_2}$ by assigning the corresponding weights and taking a sum over $j_{12}$ of a matrix element $A$, as prescribed in \cite{KR}. In this case, summation over the gleam colored by $j_{12}$ is done in accordance to the tensor product decomposition of $V_{j_1}\ot V_{j_2}$
\begin{equation}\label{vj}
V_{j_1}\ot V_{j_2}=\bigoplus\limits_{j_{12}=|j_1-j_2|}^{j_1+j_2}V_{j_{12}}.
\end{equation}
\par

The relationship between the two diagrams can be intuitively understood through the lens of the bulk-boundary correspondence. In string diagram calculus, which is essentially an edge-state model, the degrees of freedom are excitations localized at the ends of the strings and propagating along the strands. Their interactions are determined by intertwiners and braiding, depending on the specific scattering processes. In contrast, the shadow world graphical calculus is a face-state model, where the degrees of freedom reside in the bulk, and interactions between faces are governed by fusion rules arising from tensor product decomposition. The bulk-boundary correspondence asserts that edge excitation scattering data can be recoupled into constraints on the surrounding faces. Or, in other words, scattering data for the edge-model and fusion data for the face-model provide two dual descriptions of the same underlying representation-theoretic structure.

\begin{figure}[h]
\begin{center}
	{\includegraphics[width=420pt]{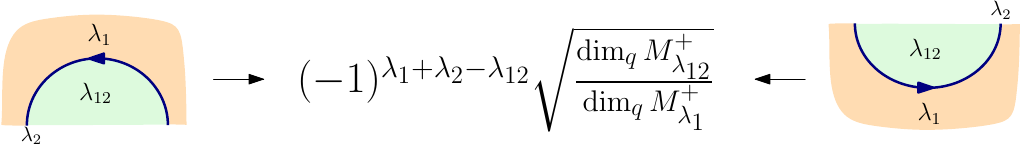}}
	\caption{Graphical representation of cups and caps in the shadow world.}
	\label{SWcc}
\end{center}
\end{figure}

Now, let us apply this principle to Verma modules. Taking a trace over a strand colored by a finite-dimensional representation contributes the same factor as in a trace over a strand colored by a Verma module, the only difference being the weight space over which the summation is being performed. Hence, we can assume that cups and caps in the shadow world should be assigned weights similar to the finite-dimensional case, but with proper quantum dimensions, as depicted in Figure \ref{SWcc}.
\par
Here
$$
\mathrm{dim}_q\,M^+_\ld=\sum\limits_{\substack{a=\ld-k \\ 
k\in\mathbb{Z}_{\geq 0}}} q^{2a}=\frac{q^{2\ld+1}}{q-q^{-1}},
$$
which is analogous to formula (\ref{vj}), except that the summation is taken over the weight space of $M^+_{\ld}$ rather than $V_j$. Once taking a trace of a braid colored by Verma modules of the highest weight, the fusion data is determined by the tensor product decomposition of $M^+_{\ld_1}\ot M^+_{\ld_2}$, with the summation over $\ld_{12}$ carried out according to (\ref{decomp}). This reflects the underlying assumption that the duality between the shadow world model and the $\mathcal{R}$-matrix model extends beyond the finite-dimensional case, capturing a broader class of representations.\par

The assignment shown in Figure \ref{SWcc} shares the same limitations as that for string diagrams. It can be applied only to traces of braids, while other topological moves remain unavailable, i.e. rotating a crossing that described by a corresponding $q6j$-symbol. Furthermore, as discussed in Section \ref{chev}, it is more natural to color gleams by polarized modules. It is straightforward to introduce gleams of two types, one for Verma modules of the highest weight and one for Verma modules of the lowest weight, and then infer that a similar argument holds, but one should use quantum dimensions of polarized modules instead. However, when taking the trace of a braid colored by polarized modules, the fusion data must be reconstructed from the existing results for the GM series.

\subsection{The functional $\mathsf{SW}_\infty$}
\subsubsection{Definition of $\mathsf{SW}_\infty$}
Let us define a functional $\mathsf{SW}_\infty$, which is a model based on the shadow world graphical calculus considered in Section \ref{shadworld}. 
\begin{definition}\label{defsw}
Consider the diagram of a braid representation of a knot $K$, which is a projection of $K$ onto $\mathbb{R}^2$ with crossings being taken into account. We color this diagram by the highest weight $\ld$ of a Verma module $M^+_\ld$. The complement to the diagram of a knot $K$ in the plane $\mathbb{R}^2$ is given by gleams, which we color by highest weights of Verma modules in accordance to the following inductive rule:
\begin{itemize}
\item The external gleam is colored by $\ld=0$.
\item Those gleams, which can be reached from the external gleam by crossing a single string, we associate the color of the crossed string.
\item Any other internal gleam, which is reached from a gleam colored by $\ld_1$ by crossing a string of color $\ld_2$, is colored by $\ld_{12}=\ld_1+\ld_2-J$, where $J\in\mathbb{Z}_{\geq_0}$. For $n$ internal gleams of this type, there will be different numbers $J_k$, $k=1,\ldots,n$.
\end{itemize}
Let us denote such a coloring of a plane $\mathbb{R}^2$ as $D_{(K,M^+_\ld)}^{\overline{J}}$, where $\overline{J}=(J_1,J_2,\ldots,J_n)$ is called the state. Functional $\mathsf{SW}_\infty[K,M^+_\ld]$ is defined as follows
\begin{itemize}
\item It takes a knot $K$ and a representation $M^+_\ld$ as an input and assigns a coloring $D_{(K,M^+_\ld)}^{\overline{J}}$.
\item To each elementary graph in a coloring $D_{(K,M^+_\ld)}^{\overline{J}}$ it assigns weights in accordance to the rules given in Figures \ref{crossing3}, \ref{trivv}, \ref{SWcc}.
\item The obtained expression should be summed over all states $\overline{J}$, where for each $J_k$ the summation should be done in accordance with the tensor product decomposition
\begin{equation}\label{tpr}
M^+_{\ld_1}\ot M^+_{\ld_2}=\bigoplus\limits_{J=0}^\infty M^+_{\ld_1+\ld_2-J}.
\end{equation}
\end{itemize}
\end{definition}
By definition, this functional is invariant with respect to transformations described in Sections \ref{ortsect}, \ref{racahsect}, \ref{BEsect}, and \ref{YBsect}. The class of knots to which the above functional applies is not yet known, but to reverse-engineer the fusion data, we follow \cite{P1} and below consider an example of a positive braid knot.\par

It is possible to apply this functional to links and knotted graphs with trivalent vertices, where diagrams are colored by Verma modules of the highest weight and finite-dimensional representations, with (\ref{tpr}) modified correspondingly. For knotted graphs with multiple components and links, the colorings of the components of the diagram should be chosen generically to avoid possible complications arising from the tensor product decomposition.

\subsubsection{Example: the trefoil}
\begin{figure}[h]
\begin{center}
	{\includegraphics[width=400pt]{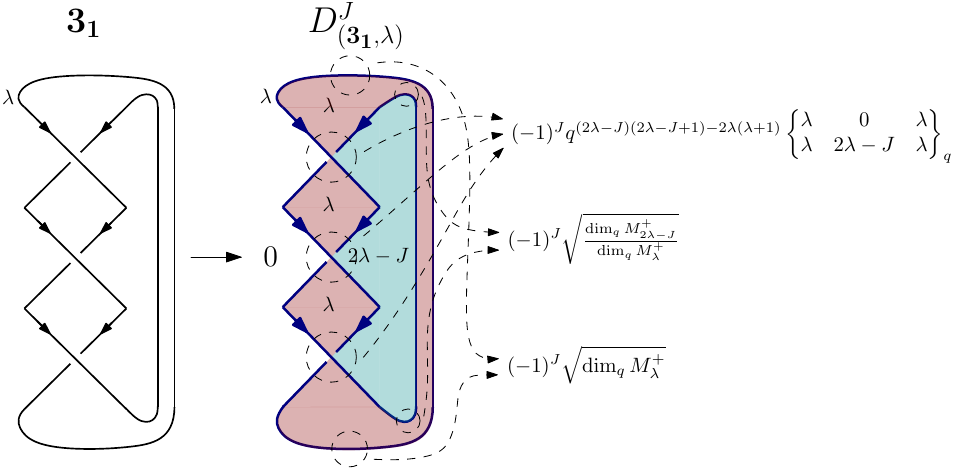}}
	\caption{Trefoil knot $\mathbf{3_1}$, its diagram and coloring $D^J_{(\mathbf{3_1},M^+_\ld)}$. In the right hand side we assign weights to each elementary graph in the considered coloring.}
	\label{trefoilqq}
\end{center}
\end{figure}

Let us compute $\mathsf{SW}_\infty[\mathbf{3_1},M^+_\ld](x,q)$. We take a positive braid representation for $\mathbf{3_1}$ and assign a coloring $D^J_{(\mathbf{3_1},M^+_\ld)}$, as depicted in Figure \ref{trefoilqq}.
\par
As a result, we obtain
\begin{eqnarray}\label{swp}
&\mathsf{SW}_\infty[\mathbf{3_1},M^+_\ld](x,q)=\sum\limits_{J=0}^\infty \left((-1)^{J}q^{(2\ld-J)(2\ld-J+1)-2\ld(\ld+1)}\qs{\ld}{0}{\ld}{\ld}{2\ld-J}{\ld}\right)^3 \times\\
&\times\left((-1)^{J}\sqrt{\frac{\mathrm{dim}_q \,M^+_{2\ld-J}}{\mathrm{dim}_q\, M^+_{\ld}}} \right)^2 \left((-1)^{J}\sqrt{\mathrm{dim}_q \,M^+_{\ld}} \right)^2=\frac{q^{6\ld^2+4\ld+1}}{q-q^{-1}}\sum\limits_{J=0}^\infty (-1)^J q^{3J^2-5J-12J\ld},\nonumber
\end{eqnarray}
where we used formulae presented in Section \ref{sectionpropq6j}.

\subsubsection{Recovering $F_{\mathbf{3_1}}(x,q)$}
Recall that in \cite{GM} the following result was obtained
\begin{equation}\label{fktref}
    F_{\mathbf{3_1}}(x,q)=\frac{q}{2}\sum\limits_{m=1}^\infty \epsilon_m(x^{\frac{m}{2}}-x^{-\frac{m}{2}})q^{\frac{m^2-1}{24}},
\end{equation}
where
\begin{equation}\label{epstref}
\epsilon_m=
\begin{cases}
-1, & \text{if $m=1$ or $11$ $(\mrm{mod}\, 12)$},\\
+1, & \text{if $m=5$ or $7$ $(\mrm{mod}\, 12)$},\\
0, & \text{otherwise}
\end{cases}.
\end{equation}
Let us focus on the negative expansion $F^-_{K}(x,q)$, where
\begin{equation}\label{fullfk}
F_K(x,q)=\frac{1}{2}(F^-_{K}(x,q)-F^-_{K}(x^{-1},q)).
\end{equation}

Consider (\ref{swp}), substitute $q\to q^{1/2}$ and denote $x:=q^{2\ld+1}$. We get
\begin{equation}\label{swtref}
\mathsf{SW}_\infty[\mathbf{3_1},M^+_\ld](x,q)=\frac{q^{3C_\ld+1}x^{-\frac{1}{2}}}{q^{\frac{1}{2}}-q^{-\frac{1}{2}}}\sum\limits_{J=0}^\infty (-1)^J q^{\frac{3J^2+J}{2}}x^{-3J}.
\end{equation}

Up to a normalization factor, series obtained in (\ref{swtref}) constitutes only half of the terms in the negative expansion $F^-_{\mathbf{3_1}}(x,q)$ of (\ref{fktref}). From Section \ref{chev} it is evident that one needs to color the knot using a polarized module $P_\ld$ of $\widetilde{\Uq}$. When taking a trace of a braid, repetitive junction of cups and caps is linear in quantum dimensions, therefore, $\mathsf{SW}_\infty[\mathbf{3_1},P_\ld]$ should contain $\mathsf{SW}_\infty[\mathbf{3_1},M^+_{\ld}]+\mathsf{SW}_\infty[\mathbf{3_1},M^-_{-\ld}]$, which also makes sense considering the tensor product decomposition 
$$
P_{\ld_1}\otimes P_{\ld_2}\supset \bigoplus\limits_{J=0}^\infty P_{\ld_1+\ld_2-J}.
$$
Let us compute $\mathsf{SW}_\infty[\mathbf{3_1},M^-_{-\ld}](x,q)$. Note that in Figure \ref{trefoilqq} weights of crossings depend on $q6j$-symbols and classical Casimir elements $C_\ld=\ld(\ld+1)$. Since $P_\ld$ is an incarnation of a large color limit of a finite-dimensional representation, these weights should remain the same for the "lower" part of $P_\ld$, as they would have remained for the "lower" part of a finite-dimensional representation. However, quantum dimensions $\mathrm{dim}_q\,M^+_{\ld}$ should be substituted by $\mathrm{dim}_q\,M^-_{-\ld}$, where
$$
\mathrm{dim}_q\,M^-_{-\ld}=-\frac{q^{-2\ld-1}}{q-q^{-1}},
$$
which is a meromorphic continuation of $\mrm{Tr}_{M^-_{-\ld}}(K)$ to a unit disk $|q|<1$. In order to apply the functional $\mathsf{SW}_\infty[\bullet,M^-_{-\ld}]$ one should also take into account that
$$
M^-_{-\ld_1}\ot M^-_{-\ld_2}=\bigoplus\limits_{J=0}^\infty M^-_{-\ld_1-\ld_2+J}.
$$
We substitute $q\to q^{1/2}$, $x=q^{2\ld+1}$ and obtain
\begin{equation}\label{swtrefm}
\mathsf{SW}_\infty[\mathbf{3_1},M^-_{-\ld}](x,q)=-\frac{q^{3C_\ld+2}x^{-\frac{5}{2}}}{q^{\frac{1}{2}}-q^{-\frac{1}{2}}}\sum\limits_{J=0}^\infty (-1)^J q^{\frac{3J^2+5J}{2}}x^{-3J}.
\end{equation}
Now we see that the following identity holds 
$$
F^-_{\mathbf{3_1}}(x,q)=q^{-C_\ld \omega(\mathbf{3_1})}(q^{\frac{1}{2}}-q^{-\frac{1}{2}})\left(\mathsf{SW}_\infty[\mathbf{3_1},M^+_{\ld}]+\mathsf{SW}_\infty[\mathbf{3_1},M^-_{-\ld}]\right)(x,q),
$$
where $\omega(K)$ is the writhe of a knot $K$.\par

\subsection{The category of polarized modules}
\subsubsection{The conjecture}
In \cite{RT2}, Reshetikhin and Turaev constructed a non-perturbative completion of the Chern-Simons partition function, which is a surgery TQFT that "upgrades" a link invariant given by the Reshetikhin-Turaev functor to a $3$-manifold invariant, also known as the Witten-Reshetikhin-Turaev invariant. By analogy, Gukov-Putrov-Pei-Vafa series $\hat{Z}_a(M_3;q)$ are expected to stand in the same relation to Gukov-Manolescu series $F_K(x,q)$. As a step toward a non-perturbative definition of $\hat{Z}_a(M_3;q)$, it is natural to first seek a functorial formulation of $F_K(x,q)$, starting with the identification of a suitable target category for the Reshetikhin-Turaev functor. Such a formulation would also be a progress toward a categorification of $F_K(x,q)$, in the same spirit that Khovanov homology provides a categorification of the Jones polynomial \cite{Kh}.\par

\begin{figure}[h]
\begin{center}
	{\includegraphics[width=280pt]{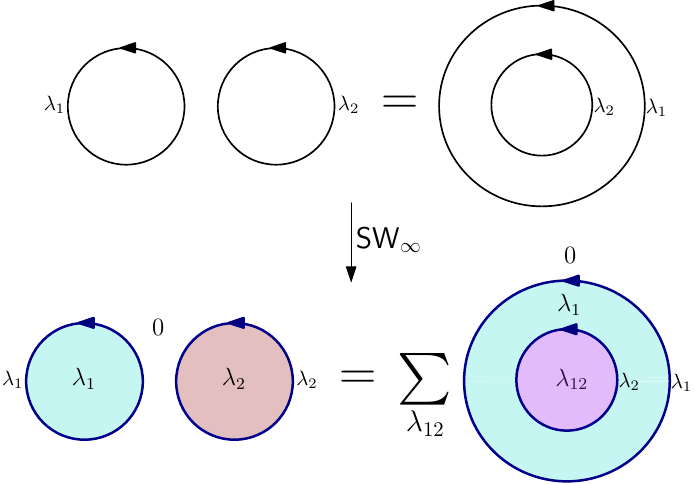}}
	\caption{Two unknots colored by $M^+_{\ld_1}$ and $M^+_{\ld_2}$, where we apply the functional $\mathsf{SW}_\infty$ in accordance to Definition \ref{defsw}: we produce a coloring of a plane and then sum over all states.}
	\label{func}
\end{center}
\end{figure}

Consider the identity depicted in Figure \ref{func}, where we color two unknots by $M^+_{\ld_1}$ and $M^+_{\ld_2}$.
\par
After applying the functional $\mathsf{SW}_\infty$ this translates into
$$
\mathrm{dim}_q\,M^+_{\ld_1}\mathrm{dim}_q\,M^+_{\ld_2}=\sum\limits_{J=0}^{\infty}\mathrm{dim}_q\,M^+_{\ld_1+\ld_2-J},
$$
which indeed holds, meaning that this functional respects the tensor product rules (\ref{tpr}) and
$$
\mathsf{SW}_\infty[K_1,M^+_{\ld_1}]\mathsf{SW}_\infty[K_2,M^+_{\ld_2}]=\mathsf{SW}_\infty[K_1\sqcup K_2,M^+_{\ld_1}\otimes M^+_{\ld_2}].
$$
In other words, if we were to consider the shadow world as a functor, which acts from the category of knots, subject to a certain limited set of transformations, to a braided monoidal version of the $q$-analog of BGG category $\mathcal{O}$, which includes Verma modules of the highest weight and concerns only generic combinations of highest weights, then the functional $\mathsf{SW}_\infty$ would be functorial. However, because of the aforementioned limitations, such a functional would fail to define a genuine topological invariant. For this reason, in the present work we refrain from formulating the shadow world functorially in terms of Verma modules of the highest weight alone.\par

Now, consider the identity in Figure \ref{func}, where we color two unknots by $P_{\ld_1}$ and $P_{\ld_2}$. As has been mentioned previously, $\mathsf{SW}_\infty[K,P_\ld]$ should contain $\mathsf{SW}_\infty[K,M^+_{\ld}]+\mathsf{SW}_\infty[K,M^-_{-\ld}]$, which, again, makes sense considering that
$$
P_{\ld_1}\otimes P_{\ld_2}\supset \bigoplus\limits_{J=0}^\infty P_{\ld_1+\ld_2-J},
$$
but $P_{\ld_1}\otimes P_{\ld_2}$ also contains summands with submodules coming from tensor products of Verma modules of the highest weight with Verma modules of the lowest weight. Hence, the Definition \ref{defsw} of the functional $\mathsf{SW}_\infty$ should be modified accordingly, where more states should be included in the summation.
\begin{conjecture}\label{conjecture}
There exists a quotient category of $\widetilde{\Uq}$-modules, Grothendieck ring of which has relations
$$
[P_{\ld_1}]\stackrel{1}{\cdot}[P_{\ld_2}]=\sum\limits_{J=0}^\infty[P_{\ld_1+\ld_2-J}]+\sum\limits_{J=0}^\infty[P_{\ld_1-\ld_2+J}],
$$
$$
[P_{\ld_1}]\stackrel{2}{\cdot}[P_{\ld_2}]=\sum\limits_{J=0}^\infty[P_{\ld_1+\ld_2-J}]+\sum\limits_{J=0}^\infty[P_{\ld_2-\ld_1+J}].
$$
\end{conjecture}
This quotient category is no longer a representation category, but should be understood in the sense of \cite{AP}, where the authors study categories endowed with a "reduced" tensor product structure arising from the character formulas of the modules in question. This is analogous to the construction of the Witten-Reshetikhin-Turaev invariants, where the category of representations of the small quantum group $u_q(\mathfrak{sl}_2)$ is semisimplified as a consequence of the indecomposable modules having vanishing quantum dimensions.\par

The definition of $\mathsf{SW}_\infty[\bullet,P_\ld]$ is analogous to Definition \ref{defsw}, except that (\ref{tpr}) should be modified in accordance to the Conjecture \ref{conjecture}. The identity in Figure \ref{func} for two unknots colored by $P_{\ld_1}$ and $P_{\ld_2}$ translates into 
\begin{eqnarray*}
\mathrm{dim}_q\,P_{\ld_1}\mathrm{dim}_q\,P_{\ld_2}&=&\sum\limits_{J=0}^{\infty}\mathrm{dim}_q\,P_{\ld_1+\ld_2-J}+\sum\limits_{J=0}^{\infty}\mathrm{dim}_q\,P_{\ld_1-\ld_2+J}\\
&=&\sum\limits_{J=0}^{\infty}\mathrm{dim}_q\,P_{\ld_1+\ld_2-J}+\sum\limits_{J=0}^{\infty}\mathrm{dim}_q\,P_{\ld_2-\ld_1+J},
\end{eqnarray*}
where
$$
\mathrm{dim}_q\,P_{\ld}=\frac{q^{2\ld+1}-q^{-(2\ld+1)}}{q-q^{-1}},
$$
which holds as a meromorphic continuation of the corresponding geometric series to a unit disk $|q|<1$.\par

\subsubsection{Reexamining the example: the trefoil}
Let us compute $\mathsf{SW}_\infty[\mathbf{3_1},P_\ld](x,q)$. Note, that there is no ambiguity in the definition of $\mathsf{SW}_\infty[\bullet,P_\ld]$, since $[P_{\ld}]\stackrel{1}{\cdot}[P_{\ld}]=[P_{\ld}]\stackrel{2}{\cdot}[P_{\ld}]$.  After the proper substitution of variables, we get
$$
\mathsf{SW}_\infty[\mathbf{3_1},P_\ld](x,q)=\frac{(-1)^{2\ld}q^{-C_\ld \omega(\mathbf{3_1})-\frac{1}{2}}\sum\limits_{m=1}^\infty \epsilon_m q^{\frac{m^2-1}{24}}-q^{C_\ld \omega(\mathbf{3_1})+1}\sum\limits_{m=1}^\infty \epsilon_m x^{-\frac{m}{2}}q^{\frac{m^2-1}{24}}}{q^{\frac{1}{2}}-q^{-\frac{1}{2}}},
$$
where $\epsilon_m$ is given in (\ref{epstref}). 
However, when implementing the formula (\ref{fullfk}) the additional term emerging from $\sum\limits_{J=0}^\infty[P_{J}]$ vanishes, and we get
$$
\frac{1}{2}\left(\mathsf{SW}_\infty[\mathbf{3_1},P_\ld](x,q)-\mathsf{SW}_\infty[\mathbf{3_1},P_\ld](x^{-1},q)\right)=\frac{q^{C_\ld \omega(\mathbf{3_1})}}{q^{\frac{1}{2}}-q^{-\frac{1}{2}}}F_{\mathbf{3_1}}(x,q),
$$
where $F_{\mathbf{3_1}}(x,q)$ is given by (\ref{fktref}).

\subsection{Open questions}

\begin{itemize}
\vspace{-0.5 cm}
\item Establish Conjecture \ref{conjecture}. The tensor product of a highest weight Verma module with a lowest weight Verma module is expected to decompose, at least formally, into a direct integral over \emph{dense modules} $V(\tau,\xi)$ (for definition see, e.g., \cite{Ma}). We expect that irreducible dense modules form a tensor ideal, which would allow one to consistently assign their quantum dimensions to be zero. Furthermore, a dense module containing a Verma module as a submodule is expected to have a quantum dimension coinciding with that of the Verma module, making them catergorically indistinguishable. This provides a preliminary outline of a strategy that will be investigated further in \cite{S}.

\item Since $F_K(x,q)$ is a series defined in the vicinity of $q=0$, it is tempting to speculate that the reduced tensor product structure conjectured above might be interpreted in terms of crystal bases \cite{KN,HK}. One possible approach would be to formally consider the crystal associated with the reduced tensor product of a lowest weight Verma crystal and a highest weight Verma crystal, together with the action of the commutor \cite{HeKa}. As such crystals are neither upper-normal nor lower-normal, the commutor would not act as a morphism of crystals; instead, suggesting the possibility of two distinct tensor product structures for crystals of polarized modules.

\item The graphical calculus for $q3j$- and $q6j$-symbols developed in the present paper currently offers only partial topological insight and would need to be refined into a fully functorial framework. Building on Conjecture \ref{conjecture} and the constructions in Section \ref{chev}, one could attempt a functorial reformulation of $\mathsf{SW}_\infty[\bullet,P_\ld]$. Such a reformulation may help clarify whether $\mathsf{SW}_\infty$ defines a genuine topological invariant, or instead an invariant of knot diagrams up to a restricted class of moves. 

\item A natural next step is to determine the class of knots or knot diagrams for which the functional $\mathsf{SW}_\infty$ reproduces $F_K(x,q)$ as a genuine topological invariant. In such cases, the functional could potentially provide explicit closed-form expressions for the series expansions studied in \cite{P1,P2}. Consistency checks with known results for $\hat{Z}_a(M_3;q)$
are also necessary, requiring a careful reconciliation of the definitions and normalizations used in \cite{MM,Roz} with those adopted here. 

\item The formulae for $q6j$-symbols obtained in the present paper can also serve as a starting point for constructing a Macdonald $(q,t)$-deformation of their squares \cite{AS}. Such a deformation could potentially provide a refinement of mapping class group representations arising from a conjectural TQFT expected from the properties of $\hat{Z}_a(M_3;q)$.

\item If the $q6j$-symbols for polarized modules possess the required symmetries, one may speculate that a Turaev-Viro state-sum \cite{TV} could be defined for these modules, potentially providing insight into the structure of a yet-to-be-defined surgery TQFT for $\hat{Z}_a(M_3;q)$.

\item It may be worthwhile to explore whether the approach of \cite{FSS}, which categorifies $q3j$-symbols via counts of signed isotopy classes of non-intersecting arcs in a triangle, can be adapted or extended to the context of Verma modules.

\item If we tensor Verma modules while relaxing the condition $\sum_j\ld_j\notin\frac{1}{2}\mbb{Z}_{\geq 0}$, it is plausible that the formulae obtained here can be generalized, with suitable normalizations. In this setting, the tensor product $M^+_{\ld_1} \otimes M^+_{\ld_2}$ would involve both \emph{big projective covers} and Verma modules that are no longer irreducible. For further details on such decompositions, see, e.g., \cite{Mu}.
\end{itemize}

\appendix
\renewcommand{\theequation}{A\arabic{equation}}
\setcounter{equation}{0}
\section{Quantum numbers, Gaussian binomial coefficients and identities}\label{qnumb}
For $q$ generic and $\ld\in\mbb{C}$ we define the quantum numbers
$$
\q{\ld}=\frac{q^\ld-q^{-\ld}}{q-q^{-1}}.
$$
For $n\geq 0$ we define the quantum factorial as
$$
\q{n}!=\prod\limits_{s=1}^n \frac{q^s-q^{-s}}{q-q^{-1}},\quad \q{0}!=1.
$$
Gaussian binomial coefficients are defined as
$$
\qbin{m}{n}=\frac{\q{m}!}{\q{n}!\q{m-n}!}.
$$
They satisfy the following recursion relation
\begin{equation}\label{qbinrec}
    \qbin{n}{m}=q^{\pm m}\qbin{n-1}{m}+q^{\mp(n-m)}\qbin{n-1}{m-1}.
\end{equation}
\begin{lemma}{\cite{K}}
The following identities hold
\begin{equation}\label{qbid1}
\sum\limits_{j} q^{\pm((a-b)j-(j-b))}\frac{(-1)^j}{\q{a-j}!\q{j-b}!}=(-1)^a\delta_{a,b},
\end{equation}
\begin{equation}\label{qbid2}
\sum\limits_{j}\frac{q^{\pm aj}}{\q{j}!\q{b-j}!\q{c-j}!\q{a-b-c+j}!}=q^{\pm bc}\frac{\q{a}!}{\q{b}!\q{c}!\q{a-b}!\q{a-c}!},
\end{equation}
\begin{equation}\label{qbid3}
\sum\limits_{j}q^{j(a+b-c+2)}\frac{\q{a-j}!\q{b+j}!}{\q{j}!\q{c-j}!}=q^{c(b+1)}\frac{\q{a-c}!\q{b}!\q{a+b+1}!}{\q{c}!\q{a+b-c+1}!},
\end{equation}
\begin{eqnarray}\label{qbid4}
&\sum\limits_{j}(-1)^jq^{\pm j(a-b-c+1)}\frac{\q{a-j}!}{\q{j}!\q{b-j}!\q{c-j}!}=q^{\mp bc}\frac{\q{a-b}!\q{a-c}!}{\q{b}!\q{c}!\q{a-b-c}!},\quad {a-b,a-c\geq 0}.
\end{eqnarray}
\end{lemma}
\begin{proof}
    Proof by induction, using (\ref{qbinrec}).
\end{proof}
For $a-b\in\mbb{Z}$ it will be convenient for us to denote
\begin{equation}\label{qpf}
        \q{a|b}=\begin{cases}
        \prod\limits_{j=1}^{a-b} \q{b+j}, & \text{if $a-b>0$}\\
        \frac{1}{\prod\limits_{j=1}^{b-a} \q{a+j}}, & \text{if $a-b< 0$}\\
        1, & \text{if $a=b$}
        \end{cases}
\end{equation}
It has the following properties
\begin{equation*}
    \q{a|b}=\frac{\q{a}!}{\q{b}!},\quad a,b\in{\mbb{Z}_{\geq 0}},
\end{equation*}
\begin{equation*}
    \q{a|b}=\frac{1}{\q{b|a}},
\end{equation*}
\begin{equation}\label{qpfid2}
    \q{a|b}\prod\limits_{j=0}^{c-1}\q{b-j}=\q{a|b-c},\quad c\geq 0,
\end{equation}
\begin{equation}\label{qpfid3}
    \frac{\q{a|b}}{\prod\limits_{j=1}^{c-b}\q{b+j}}=\q{a|c},\quad c-b\geq 0
\end{equation}

\section{Lemmas implementing the Racah method}\label{appp}
Now we will present some technical lemmas used throughout derivation of the $q6j$-symbols. In what follows we consider values of the parameters in accordance to their values in definition (\ref{qsdef}). More precisely, $\ld_1,\ld_2,\ld_3\in\mbb{C}$, they and their sums are generic. Also, $\ld_2+\ld_3-\ld_{23},\ld_1+\ld_2-\ld_{12}, \ld_1+\ld_2+\ld_3-\ld\in\mbb{Z}_{\geq 0}$, and $\ld_1+\ld_{23}-\ld\geq 0$, $\ld_{12}+\ld_3-\ld\geq 0$. The following lemmas hold.

\begin{lemma}[Step 2]\label{lemstep2}
For $\ld_3-s,z\in\mbb{Z}_{\geq 0}$
\begin{eqnarray*}
&\sum\limits_{p\in\mbb{Z}}\frac{q^{p(\ld_3+\ld_2+\ld_1+2\ld_{23}-\ld_{12}-z-s+2)}}{\q{p}!\q{\ld_{12}-\ld_1-\ld_{23}+z+s-p}! \prod\limits_{j=0}^{\ld_3-s+p-1}\q{2\ld_3-j}\prod\limits_{j=0}^{\ld_2-\ld_{23}+s-p-1}\q{2\ld_2-j}}=\\
&=q^{(\ld_2+\ld_{23}-s+1)(\ld_{12}-\ld_1-\ld_{23}+z+s)}\frac{\q{\ld_2+\ld_3+\ld_{23}+1|\ld_1+\ld_2+\ld_3+2\ld_{23}-\ld_{12}-z-s+1}}{\q{\ld_{12}-\ld_1-\ld_{23}+z+s}!\prod\limits_{j=0}^{\ld_3-\ld_1-\ld_{23}+\ld_{12}+z-1}\q{2\ld_3-j}\prod\limits_{j=0}^{\ld_2-\ld_{23}+s-1}\q{2\ld_2-j}}.
\end{eqnarray*}
\end{lemma}
\begin{proof}
First consider identity (\ref{qbid3}) for finite dimensional irreducible representations, where
\begin{equation*}
    a\to \ld_3+s,\quad b\to \ld_2+\ld_{23}-s,\quad c\to\ld_{12}-\ld_1-\ld_{23}+z+s,
\end{equation*}
which gives
\begin{eqnarray*}
&\sum\limits_{p\in\mbb{Z}}q^{p(\ld_2+\ld_3+\ld_1+2\ld_{23}-\ld_{12}-z-s+2)}\frac{\q{\ld_3+s-p}!\q{\ld_2+\ld_{23}-s+p}!}{\q{p}!\q{\ld_{12}-\ld_1-\ld_{23}+z+s-p}!}=\\
&=q^{(\ld_2+\ld_{23}-s+1)(\ld_{12}-\ld_1-\ld_{23}+z+s)}\frac{\q{\ld_3+\ld_2+\ld_{23}+1}!\q{\ld_2+\ld_{23}-s}!\q{\ld_3+\ld_1+\ld_{23}-\ld_{12}-z}!}{\q{\ld_{12}-\ld_1-\ld_{23}+z+s}!\q{\ld_3+\ld_2+\ld_1+2\ld_{23}-\ld_{12}-z-s+1}!}.
\end{eqnarray*}
Dividing both sides by $\q{2\ld_2}!\q{2\ld_3}!$ gives the result, since this identity admits analytical continuation to generic $\ld_1,\ld_2,\ld_3\in\mbb{C}$.
\end{proof}

\begin{lemma}[Step 3]\label{lemstep3}
For $\ld_3-s,z\in\mbb{Z}_{\geq 0}$
\begin{eqnarray*}
    &\sum\limits_{p\in\mbb{Z}}(-1)^{p} q^{p(\ld_3+\ld_2+\ld_{23}-\ld-\ld_{12}-z-s-1)}\frac{\q{\ld_2-\ld_1+\ld_{12}|\ld_2+\ld_{23}-s-z-p}}{\q{\ld_3-s-p}!\q{p}!\prod\limits_{j=0}^{\ld_{12}-\ld+s+p-1}\q{2\ld_{12}-j}}=\\
    &=q^{(\ld_3-s)(\ld_2+\ld_{23}-z-s)}\frac{\q{\ld_2-\ld_1+\ld_{12}|\ld_2+\ld_{23}-s-z}\prod\limits_{j=0}^{\ld_3-s-1}\q{\ld+\ld_{12}-\ld_2-\ld_{23}+z-j}}{\q{\ld_3-s}!\prod\limits_{j=0}^{\ld_{12}+\ld_3-\ld-1}\q{2\ld_{12}-j}}.
\end{eqnarray*}
\end{lemma}
\begin{proof}
First consider identity (\ref{qbid4}) for finite dimensional irreducible representations, where
\begin{equation*}
    a\to \ld+\ld_{12}-s,\quad b\to \ld_3-s,\quad c\to \ld_2+\ld_{23}-s-z,
\end{equation*}
which gives
\begin{eqnarray*}
    &\sum\limits_{p\in\mbb{Z}}(-1)^{p} q^{-p(\ld+\ld_{12}-\ld_2-\ld_3-\ld_{23}+z+s+1)}\frac{\q{\ld+\ld_{12}-s-p}!}{\q{p}!\q{\ld_3-s-p}!\q{\ld_2+\ld_{23}-s-z-p}!}=\\
    &=q^{(\ld_3-s)(\ld_2+\ld_{23}-z-s)}\frac{\q{\ld-\ld_3+\ld_{12}}!\q{\ld+\ld_{12}-\ld_2-\ld_{23}+z}!}{\q{\ld_3-s}!\q{\ld_2+\ld_{23}-s-z}!\q{\ld+\ld_{12}-\ld_2-\ld_3-\ld_{23}+s+z}!}.
\end{eqnarray*}
Multiplying both sides by $\frac{\q{\ld_2-\ld_1+\ld_{12}}!}{\q{2\ld_{12}}!}$ gives the result, since this identity admits analytical continuation to generic $\ld_1,\ld_2,\ld_3\in\mbb{C}$.
\end{proof}

\begin{lemma}[Step 5]\label{lemstep5}
For $r\in\mbb{Z}_{\geq 0}$
\begin{eqnarray*}
   &\sum\limits_{p\in\mbb{Z}}(-1)^{p}\frac{q^{p(\ld_3-\ld_2-\ld_1+\ld+r+1)}}{\q{p}!\q{\ld_1+\ld_2-\ld_{12}-r-p}!\q{\ld_1+\ld_{23}-\ld-r-p}!\prod\limits_{j=0}^{\ld_3-\ld_1-\ld_{23}+\ld_{12}+r+p-1}\q{2\ld_3-j}}=\\
    &=q^{-(\ld_2+\ld_1-\ld_{12}-r)(\ld_1+\ld_{23}-\ld-r)}\frac{\prod\limits_{j=0}^{\ld_1+\ld_2-\ld_{12}-r-1}\q{\ld_3+\ld-\ld_{12}-j}}{\q{\ld_1+\ld_2-\ld_{12}-r}!\q{\ld_1+\ld_{23}-\ld-r}!\prod\limits_{j=0}^{\ld_2+\ld_3-\ld_{23}-1}\q{2\ld_3-j}}.
\end{eqnarray*}
\end{lemma}
\begin{proof}
Proof is similar to the one of Lemma \ref{lemstep3}. We use (\ref{qbid4}) for
\begin{equation*}
    a\to \ld_3+\ld_1+\ld_{23}-\ld_{12}-r,\quad b\to \ld_1+\ld_2-\ld_{12}-r,\quad c\to \ld_1+\ld_{23}-\ld-r,
\end{equation*}
and divide both sides by $\q{2\ld_3}!$.
\end{proof}

\begin{lemma}[Step 6]\label{lemstep6}
For $\ld_3-u,r\in\mbb{Z}_{\geq 0}$,
\begin{eqnarray*}
    &\sum\limits_{p\in\mbb{Z}}(-1)^p\frac{q^{-p(\ld_2-\ld_3+\ld_{23}-r+1)}}{\q{r-p}!\q{p}!\q{\ld_3-u+r-p}!\prod\limits_{j=0}^{\ld_2-\ld_{23}+u-r+p-1}\q{2\ld_2-j}}=\\
    &=q^{r(\ld_3-u+r)}\frac{\prod\limits_{j=0}^{r-1}\q{\ld_2+\ld_{23}-\ld_3-j}}{\q{r}!\q{\ld_3-u+r}!\prod\limits_{j=0}^{\ld_2-\ld_{23}+u-1}\q{2\ld_2-j}}.
\end{eqnarray*}
\end{lemma}
\begin{proof}
Proof is similar to the one of Lemma \ref{lemstep3}. We use (\ref{qbid4}) for
\begin{equation*}
    a\to \ld_2+\ld_{23}-u+r,\quad b\to r,\quad c\to \ld_3-u+r,
\end{equation*}
and divide both sides by $\q{2\ld_2}!$.
\end{proof}

\begin{lemma}[Step 7]\label{lemstep7}
For $r\in\mbb{Z}_{\geq 0}$
\begin{eqnarray*}
    &\sum\limits_{p\in\mbb{Z}}q^{p(\ld_3+\ld+\ld_{12}+r+1)}\frac{\q{\ld_2+\ld_3+\ld_{23}+1|\ld_1+\ld_2+2\ld_{23}-\ld_{12}+p-r+1}\q{\ld_1-\ld_2+\ld_{12}|\ld+\ld_{12}-\ld_2-\ld_{23}-p+r}}{\q{p}!\q{\ld_{12}-\ld_1-\ld_{23}+\ld_3+p-r}!}=\\
    &=q^{(\ld_{12}-\ld_1-\ld_{23}+\ld_3+r)(\ld+\ld_{12}-\ld_2-\ld_{23}+r)}\frac{\q{\ld_1-\ld_2+\ld_{12}|\ld+\ld_{12}-\ld_2-\ld_{23}+r}\prod\limits_{j=0}^{\ld_{12}-\ld_{23}+\ld_3-\ld_1+r-1}\q{\ld_3+\ld+\ld_{12}+r+1-j}}{\q{\ld_{12}-\ld_1-\ld_{23}+\ld_3+r}!}.
\end{eqnarray*}
\end{lemma}
\begin{proof}
Consider identity (\ref{qbid2}) for finite dimensional irreducible representations, where
\begin{equation*}
    a\to \ld_3+\ld+\ld_{12}+r+1,\quad b\to \ld_{12}-\ld_1-\ld_{23}+\ld_3+r,\quad c\to \ld+\ld_{12}-\ld_2-\ld_{23}+r,
\end{equation*}
which gives
\begin{eqnarray*}
    &\sum\limits_{p\in\mbb{Z}}\frac{q^{p(\ld_3+\ld+\ld_{12}+r+1)}}{\q{p}!\q{\ld_{12}-\ld_1+\ld_3-\ld_{23}-p+r}!\q{\ld+\ld_{12}-\ld_2-\ld_{23}-p+r}!\q{\ld_2+\ld_1+2\ld_{23}-\ld_{12}+1+p-r}!}=\\
    &=q^{(\ld_{12}-\ld_1-\ld_{23}+\ld_3+r)(\ld+\ld_{12}-\ld_2-\ld_{23}+r)}\frac{\q{\ld_3+\ld+\ld_{12}+r+1}!}{\q{\ld_{12}-\ld_1-\ld_{23}+\ld_3+r}!\q{\ld+\ld_{12}-\ld_2-\ld_{23}+r}!\q{\ld+\ld_{23}+\ld_1+1}!\q{\ld_3+\ld_2+\ld_{23}+1}!}.
\end{eqnarray*}
Multiplying both sides by $\frac{\q{\ld_1-\ld_2+\ld_{12}}!}{\q{\ld_3+\ld_2+\ld_{23}+1}!}$ gives the result, since this identity admits analytical continuation to generic $\ld_1,\ld_2,\ld_3\in\mbb{C}$. Note that the presence of $\q{\ld_{12}-\ld_1-\ld_{23}+\ld_3+r}!$ implies
\begin{equation*}
    \q{\ld_3+\ld+\ld_{12}+r+1|\ld+\ld_{23}+\ld_1+1}=\prod\limits_{j=0}^{\ld_{12}-\ld_{23}+\ld_3-\ld_1+r-1}\q{\ld_3+\ld+\ld_{12}+r+1-j}.
\end{equation*}
\end{proof}

\section{$q$-Pochhammer symbols, $q$-hypergeometric series and related identities}
We will be following monography \cite{GR} and references therein. The $q$-Pochhammer symbol is defined as
$$
(a;q)_n=
\begin{cases}
1,& n=0,\\
\prod\limits_{j=0}^{n-1}(1-aq^j), & n=1,2,\ldots,
\end{cases}
$$
and
$$
(a;q)_\infty=\prod\limits_{j=0}^{\infty}(1-aq^j).
$$
They satisfy the following identities
\begin{equation}\label{id1}
(a;q)_n=\frac{(a;q)_\infty}{(aq^n;q)_\infty},
\end{equation}
\begin{equation}\label{id2}
    (a;q)_{n-k}=(-qa^{-1})^kq^{\frac{k(k-1)}{2}-nk}\frac{(a;q)_n}{(a^{-1}q^{1-n};q)_k},
\end{equation}

We will need the following identities, relating $q$-Pochhammer symbols to quantum numbers
\begin{equation}\label{qid1}
\q{k}!=\frac{(-1)^k q^{-\frac{k}{2}(k+1)}}{(q-q^{-1})^k}(q^2;q^2)_k,
\end{equation}
\begin{equation}\label{qid2}
\p{k-1}{J-j}=\frac{q^{\frac{k}{2}(2J+1-k)}}{(q-q^{-1})^k}(q^{-2J};q^2)_k.
\end{equation}

We will use notation
$$
(a_1,\dots,a_k;q)_n=(a_1;q)_n\dots(a_k;q)_n.
$$
We define $q$-hypergeometric series as
$$
{}_r \phi_s\left[\begin{matrix}
a_1,& \dots, & a_r\\
b_1,& \dots, & b_s
\end{matrix};q,z\right]=\sum\limits_{n=0}^\infty\frac{(a_1,\dots,a_r;q)_n}{(q,b_1,\dots,b_s;q)_n}\left((-1)^n q^{\frac{n(n-1)}{2}} \right)^{1+s-r}z^n.
$$
The following identity holds
\begin{equation}\label{hid1}
{}_2 \phi_1\left[\begin{matrix}
a, b\\
c
\end{matrix};q,\frac{c}{ab}\right]=\frac{(\frac{c}{a}, \frac{c}{b};q)_\infty}{(c,\frac{c}{ab};q)_\infty}.
\end{equation}

\section{Proof of Lemma \ref{shadowtransit}}\label{prooftr}
Identity (\ref{stid1}) was obtained in (\ref{step02}). Identity (\ref{stid2}) is similar to (\ref{step01}). Let us briefly sketch the proof for the identity (\ref{stid3}). One starts with its left hand side. Use identity (\ref{rmatid1}) as
\begin{eqnarray*}
    &\qt{\ld_{13}}{\ld_2}{\ld_{123}}{a_{13}}{a_2}{a_{123}}^\psi=(-1)^{\ld_{2}+\ld_{13}-\ld_{123}}q^{C_{123}-C_2-C_{13}}\sum\limits_{a_2^\prime,a_{13}^\prime}\left(\left(\mathcal{R}^{\ld_2,\ld_{13}}\right)^{-1}\right)^{a_2^\prime,a_{13}^\prime}_{a_2,a_{13}}\qt{\ld_2}{\ld_{13}}{\ld_{123}}{a_2^\prime}{a_{13}^\prime}{a_{123}}^\psi.
\end{eqnarray*}
Then use identity (\ref{rmatid2}) to transform $$
\mathcal{R}^{\ld_3,\ld_2}\left(\mathcal{R}^{\ld_2,\ld_{13}}\right)^{-1}\to \mathcal{R}^{\ld_3,\ld_2}\left(\mathcal{R}^{\ld_2,\ld_1}\mathcal{R}^{\ld_2,\ld_3}\right)^{-1}=\left(\mathcal{R}^{\ld_2,\ld_1}\right)^{-1}
$$ 
and apply identity (\ref{stid2}) as
\begin{eqnarray*}
    \sum\limits_{a_2,a_1,a_3}\qt{\ld_2}{\ld_{13}}{\ld_{123}}{a_2}{a_{13}}{a_{123}}^\psi\qt{\ld_1}{\ld_3}{\ld_{13}}{a_1}{a_3}{a_{13}}^\psi=\sum\limits_{\ld_{12}}\qs{\ld_3}{\ld_1}{\ld_{13}}{\ld_2}{\ld_{123}}{\ld_{12}}\times\\
    \times \sum\limits^\prime_{a_2,a_1,a_3}\qt{\ld_{12}}{\ld_3}{\ld_{123}}{a_{12}}{a_3}{a_{123}}^\psi\qt{\ld_2}{\ld_1}{\ld_{12}}{a_2}{a_1}{a_{12}}^\psi.
    \end{eqnarray*}
    The last step is to apply identity (\ref{rmatid1})
    \begin{eqnarray*}
    &\sum\limits_{a_1^\prime,a_2^\prime}\left(\left(\mathcal{R}^{\ld_2,\ld_1}\right)^{-1}\right)^{a_1^\prime,a_2^\prime}_{a_1,a_2}\qt{\ld_2}{\ld_1}{\ld_{12}}{a_2^\prime}{a_1^\prime}{a_{12}}^\psi=(-1)^{\ld_{12}-\ld_1-\ld_2}q^{-C_{12}+C_1+C_2}\qt{\ld_1}{\ld_2}{\ld_{12}}{a_1}{a_2}{a_{12}}^\psi.
\end{eqnarray*}

\end{document}